\documentclass[reqno,a4paper,11pt]{article}

\usepackage{CJK,amsmath,amsthm,dsfont,amsfonts,amssymb,fancyhdr}
\usepackage{enumerate}
\usepackage{bbm,color,soul,supertabular,longtable,verbatim,extarrows}
\usepackage{titlesec}
\usepackage{graphicx}
\usepackage{cite}
\usepackage[export]{adjustbox}
\usepackage{graphicx,multirow,bm,rotating}
\usepackage[title]{appendix} 
\usepackage{booktabs,multirow,makecell}
\usepackage{tikz}
\usepackage{verbatim}
\usepackage[noblocks]{authblk}
\usepackage{ulem}
\usepackage{cancel}
\usepackage{subfigure}
\usepackage{booktabs,multirow,makecell}
\usepackage[top=2.4cm,bottom=2.2cm,left=2.6cm,right=2cm]{geometry} 
\usepackage{hyperref}
\usepackage{authblk}

\newtheorem{theorem}{Theorem}
\theoremstyle{definition}

\newtheorem{definition}[theorem]{Definition}
\theoremstyle{plain}
\newtheorem{lemma}[theorem]{Lemma}
\newtheorem{corollary}[theorem]{Corollary}
\newtheorem{remark}[theorem]{Remark}
\newtheoremstyle{noparens}
  {}{}
  {\itshape}{}
  {\bfseries}{.}
  { }
  {\thmname{#1}\thmnumber{ #2}\mdseries\thmnote{ #3}}

\theoremstyle{noparens}
\newtheorem{lemmaNoParens}[theorem]{Lemma}
\newtheorem{theoremNoParens}[theorem]{Theorem}
\newtheorem{propositionNoParens}[theorem]{Proposition}

\newtheorem*{reptheorem}{Theorem}
\title{\textbf{On the Maximum Number of Maximum Generalized $k$-Independent Sets in Forests}\thanks{This work was supported by Beijing Natural Science Foundation (Grant No. 1252010).}}

\author{
    Ning Zhang, \quad  Min Xu\thanks{\hangindent=1.5em Corresponding author. \newline {\em E-mail address:} xum@bnu.edu.cn (M. Xu).}
}

\date{}
\begin{document}
\maketitle
\vspace{-2em}

\begin{center}
    \emph{School of Mathematical Sciences, Beijing Normal University,} \\
    \emph{Key Laboratory of Mathematics and Complex Systems, Ministry of Education,} \\
    \emph{Beijing, 100875, China}
\end{center}

\renewcommand\abstractname{Abstract}
\begin{abstract}
A maximum generalized $k$-independent set in a graph is a vertex subset of maximum cardinality whose induced subgraph contains no trees of order $k$. Let $\Phi_k(G)$ denote the number of such sets in a graph $G$. Jou and Chang (2000) determined the tight upper bound of $\Phi_2(F)$ for a forest $F$ with a fixed order. Sun and Li (2023) determined the tight upper bound of $\Phi_3(F)$ for a forest $F$ under some constraints. In this paper, we investigate the general case for $k\geq2$. For a forest $F$ of order $n=pk+q$ ($n\geq1,\,0\leq q\leq k-1$), we establish a tight upper bound $\Phi_k(F)\leq k^p$. Moreover, under the condition that every component of $F$ has order at least $k$, we show that $\Phi_k(F)\leq (k-\lfloor\frac{q}{p}\rfloor)^{p-r}\cdot(k-\lceil\frac{q}{p}\rceil)^{r}$, where $0\leq r\leq p-1$ and $r\equiv q \pmod{p}$. We also characterize all extremal forests attaining these bounds. Our results extend the earlier works mentioned above.
\end{abstract}

\section{Introduction}
\quad\,\, In this paper, we consider simple, finite, and undirected graphs. Let $G=(V(G),\,E(G))$ be a graph with vertex set $V(G)$ and edge set $E(G)$. The order of $G$ is the number $n=|V(G)|$. For an integer $k\geq1$, a \textit{$k$-path} is a path of order $k$, denoted by $P_k$, and a \textit{$k$-tree} is a tree of order $k$. A \textit{forest} is a disjoint union of trees. For two graphs $G_1$ and $G_2$, we denote their disjoint union by $G_1\cup G_2$ and the disjoint union of $k$ copies of $G_1$ by $k\,G_1$. The \textit{vertex-induced subgraph} on a set $S \subseteq V(G)$ is denoted by $G[S]$. Similarly, the \textit{edge-induced subgraph} on a set $X \subseteq E(G)$ is denoted by $G[X]$. We follow \cite{1} for graph-theoretical terminology and notations not defined here.

A subset $S\subseteq V(G)$ is an \textit{independent set} if the corresponding induced subgraph $G[S]$ contains no edges. An independent set is called \textit{maximum} if it has maximum cardinality. The \textit{independence number} $\alpha(G)$ of $G$ is the cardinality of a maximum independent set of $G$. A subset $M\subseteq E(G)$ is a \textit{matching} if no two edges in $M$ share a common vertex. The \textit{matching number} $\mu(G)$ of $G$ is the cardinality of a maximum matching of $G$. A subset $C \subseteq V(G)$ is a \textit{vertex cover} if every edge of $G$ has at least one endpoint in $C$. The \textit{vertex cover number} $\tau(G)$ of $G$ is the cardinality of a minimum vertex cover of $G$. The following identity was observed by Gallai \cite{9}: $\alpha(G)+\tau(G)=|V(G)|$. Thus, $\alpha(G)+\mu(G)\leq|V(G)|$ since $\mu(G)\leq \tau(G)$. Specifically, a graph $G$ is called a \textit{K\"{o}nig–Egerv\'{a}ry graph} if $\alpha(G)+\mu(G)=|V(G)|$. The classic K\"{o}nig–Egerv\'{a}ry theorem \cite{13,14} asserts that every bipartite graph is a K\"{o}nig–Egerv\'{a}ry graph. In particular, all forests are K\"{o}nig–Egerv\'{a}ry graphs. 

Clearly, an edge is both a $2$-path and a $2$-tree. Therefore, the independent set can also be viewed as a set of vertices whose induced subgraph contains no $2$-paths or $2$-trees. In 1981, Yannakakis \cite{4} considered a vertex set whose induced subgraph contains no $3$-paths or $3$-trees, which is called a \textit{dissociation set}. The \textit{dissociation number} $diss(G)$ of $G$ is the cardinality of a maximum dissociation set of $G$. Yannakakis \cite{4} proved that computing $diss(G)$ is already NP-complete for bipartite graphs. More related results concerning the dissociation number can be found in \cite{10,11,12}. Naturally, Tu et al. \cite{2} generalized the concept of an independent set in $G$ to a $k$-\textit{independent set}, defined as a vertex set $S$ such that $G[S]$ does not contain any $k$-paths. After that, Li and Xu \cite{3} generalized the concept of an independent set in $G$ to a \textit{generalized} $k$-\textit{independent set}, defined as a vertex set $S$ such that $G[S]$ does not contain any $k$-trees. Note that $P_k$ is also a $k$-tree. Therefore, every generalized $k$-independent set is a $k$-independent set, but the converse is not true. Similarly, the $k$-\textit{independence number} $\alpha'_k(G)$ (resp. \textit{generalized} $k$-\textit{independence number} $\alpha_k(G)$) of $G$ is the cardinality of a maximum $k$-independent set (resp. maximum generalized $k$-independent set) in $G$. A $k$-\textit{matching} (resp. \textit{generalized} $k$-\textit{matching}) of $G$ is a collection of vertex-disjoint $k$-paths (resp. $k$-trees) in $G$, with the corresponding $k$-\textit{matching number} $\mu'_k(G)$ (resp. \textit{generalized} $k$-\textit{matching number} $\mu_k(G)$) denoting its maximum cardinality. A graph $G$ is a $k$-\textit{K\"{o}nig–Egerv\'{a}ry graph} if $\alpha'_k(G)+\mu'_k(G)=|V(G)|$, and it is a \textit{generalized} $k$-\textit{K\"{o}nig–Egerv\'{a}ry graph} if $\alpha_k(G)+\mu_k(G)=|V(G)|$. Specifically, any forest is both a $k$-K\"{o}nig–Egerv\'{a}ry graph and a generalized $k$-K\"{o}nig–Egerv\'{a}ry graph, as shown in \cite{2} and \cite{3}, respectively.

For $k=2$, the generalized $2$-independent set is exactly the independent set. For $k=3$, the generalized $3$-independent set is exactly the dissociation set. In approximately 1960, Erd\"{o}s and Moser posed the problems of determining the maximum number of maximum generalized $2$-independent sets among all graphs of order $n$ and of characterizing extremal graphs. This question was settled by Moon and Moser \cite{5} in 1965. Zito \cite{6} subsequently considered the same problems for the class of trees. She determined the maximum number of maximum generalized $2$-independent sets in trees of order $n$ and characterized the extremal trees that attain these bounds. Motivated by Zito's work, Tu et al. \cite{2} determined the maximum number of maximum generalized $3$-independent sets in trees of order $n$ and described all extremal trees attaining this maximum. Recently, Li and Xu \cite{3} determined the maximum number of maximum generalized $4$-independent sets in trees of order $n$ and characterized all the extremal trees. Notably, they also established some important structure theorems about maximum generalized $k$-independent sets in a tree for a general integer $k$.

It is natural to determine the maximum number of the maximum generalized $k$-independent sets in forests of order $n$ and characterize all the extremal graphs. Furthermore, it is also natural to consider this problem under the condition that each component of the forest has order at least $k$, since any maximum generalized $k$-independent set in a forest must include all vertices belonging to components with orders less than $k$. Let $\Phi_k(F)$ denote the total number of maximum generalized $k$-independent sets in a forest $F$. For $k=2$, Jou and Chang \cite{7} showed that for any forest $F$ of order $n=2p+q$ ($n\geq1,\,0\leq q\leq1$), $\Phi_2(F)\leq 2^p$, with the equality holding if and only if $F \cong p K_2 \cup q K_1$. For $k=3$, Sun and Li \cite{8} recently showed that for a forest $F$ of order $n=3p+q$ ($p\geq 2,\, 0\leq q\leq2$) with each component having order at least 3, $\Phi_3(F)\leq3^{p-q}\cdot 2^{q}$; they also showed that the equality holds if and only if $F \cong (p-q) P_3\cup q P_4$. In this paper, we consider the upper bound of $\Phi_k(F)$ for a forest $F$ for a general integer $k\geq2$, and obtain the results shown below.

\begin{theorem}\label{thm:g-zzh}
Let $k\geq 2$ be an integer$,$ and let $F$ be a forest of order $n=pk+q$ $(n\geq1,\,0\leq q\leq k-1)$. Then$,$ $\Phi_k(F)\leq k^{p}$. The equality holds if and only if $F$ consists of $p$ trees of order $k$ and an arbitrary forest of order $q$.
\end{theorem}

\begin{theorem}\label{thm:g-lsc}
Let $k\geq2$ be an integer$,$ and let $F$ be a forest of order $n=pk+q$ $(p\geq1,\,0\leq q\leq k-1),$ where each component of $F$ has order at least $k$. Then$,$ $\Phi_k(F)\leq (k-\lfloor\frac{q}{p}\rfloor)^{p-r}\cdot(k-\lceil\frac{q}{p}\rceil)^{r},$ where $0\leq r\leq p-1$ and $r\equiv q \pmod{p}$. Furthermore$,$ all the extremal graphs are described as follows$\,:$

\noindent $(1)$ If $k=2$ and $q=0,$ then the equality holds $\iff F\cong pK_2$.

\noindent $(2)$ If $k=2$ and $q=1,$ then the equality holds $\iff F$ is a forest constructed by attaching some $($possibly none$)$ of the $p-1$ copies of $K_2$ to the central vertex of a $3$-path $P$.

\noindent $(3)$ If $k\geq3$ and $p\geq q,$ then the equality holds $\iff F\cong (\bigcup\limits_{i=1}^{p-q} T_i)\cup qP_{k+1},$ where each $T_i$ is a $k$-tree.

\noindent $(4)$ If $k\geq3$ and $p<q$ with $p\neq1$ or $q\neq k-1,$ then the equality holds $\iff F\cong \bigcup\limits_{i=1}^{p} T_i,$ where each $T_i$ is a $(k+t_i)$-tree obtained from a $(k-t_i)$-path by attaching arbitrary forests $F_{i}^{1}$ and $F_{i}^{2}$ of order $t_i$ to its two ends$,$ respectively. Here$,$ $t_i=\lfloor\frac{q}{p}\rfloor$ for $1\leq i\leq p-r$ and $\lceil\frac{q}{p}\rceil$ for $p-r< i\leq p$.

\noindent $(5)$ If $k\geq3,$ $p=1$ and $q=k-1,$ then the equality holds $\iff F$ is an arbitrary $(2k-1)$-tree.
\end{theorem}

\begin{figure}[htbp]
    \centering 
    \includegraphics[width=1\textwidth]{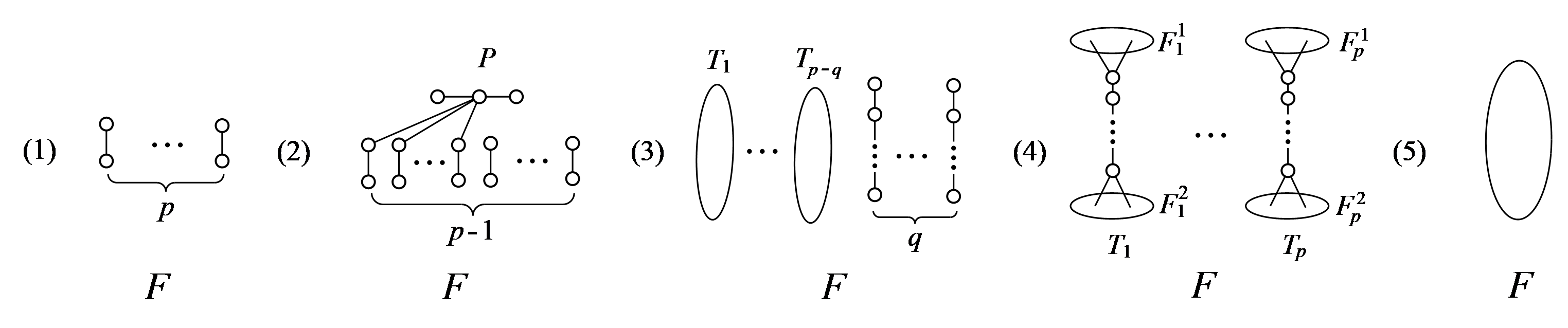}
    \caption{Extremal graphs for Theorem \ref{thm:g-lsc}.}
    \label{fig:Thm2} 
\end{figure}

This paper is organized as follows. In Section 2, we recall the structure theorems concerning the maximum generalized $k$-independent sets obtained in \cite{3} for trees and extend them to forests. In Section 3, we establish tight upper bounds on $\Phi_k(F)$ for any forest $F$ with a fixed order; these bounds are shown in Theorem \ref{thm:g-zzh} and Theorem \ref{thm:g-lsc}. Finally, in Section 4, we conclude this paper.

\section{The Structure Theorems}
\quad\,\,Throughout this paper, unless specified otherwise, let $k\geq2$ be an integer. Recall that $\alpha_k(G)$ denotes the generalized $k$-independence number of $G$, while $\mu_k(G)$ represents the generalized $k$-matching number. For any forest $F$, the relation $\alpha_k(F)+\mu_k(F)=|V(F)|$ holds. An edge $e\in E(G)$ is called $\alpha_k$-\textit{critical} if $\alpha_k(G-e)>\alpha_k(G)$ and $\mu_k$--\textit{critical} if $\mu_k(G-e)<\mu_k(G).$ Li and Xu \cite{3} established the following fundamental properties. 

\begin{lemmaNoParens}[\cite{3}]\label{lem:alpha}
If $uv$ is an $\alpha_k$-critical edge of $G,$ then every maximum generalized $k$-independent set of $G-uv$ contains both $u$ and $v,$ and $\alpha_k(G-uv)=\alpha_k(G)+1$.
\end{lemmaNoParens}

\begin{lemmaNoParens}[\cite{3}]\label{lem:mu}
If $uv$ is a $\mu_k$-critical edge of $G,$ then $uv$ is an edge of every maximum generalized $k$-matching of $G$ and $\mu_k(G-uv)=\mu_k(G)-1$.
\end{lemmaNoParens}

They also showed that an edge is $\alpha_k$-critical if and only if it is $\mu_k$-critical in a tree. Clearly, the same characterization holds for forests. Throughout the remainder of this paper, we focus exclusively on forests. Therefore, we uniformly refer to such edges as $\alpha_k$-\textit{critical edge}s for consistency.

We categorize the vertices contained in a forest $F$ into three distinct types. A vertex of $F$ is called \textit{flexible} if it is contained in some but not all maximum generalized $k$-independent sets. Let $\mathcal{F}_F=\{v:v \text{ is flexible in } F\}$. A vertex of $F$ is called \textit{static-included} if it appears in every maximum generalized $k$-independent set. Let $\mathcal{A}_F=\{v:v \text{ is static-included in } F\}$. Conversely, a vertex of $F$ is \textit{static-excluded} if it does not belong to any maximum generalized $k$-independent sets, and let $\mathcal{N}_F=\{v:v \text{ is static-excluded in } F\}$. Notably, $V(F)=\mathcal{F}_F\cup \mathcal{A}_F\cup \mathcal{N}_F$.

\begin{figure}[htbp]
    \centering 
    \includegraphics[width=0.6\textwidth]{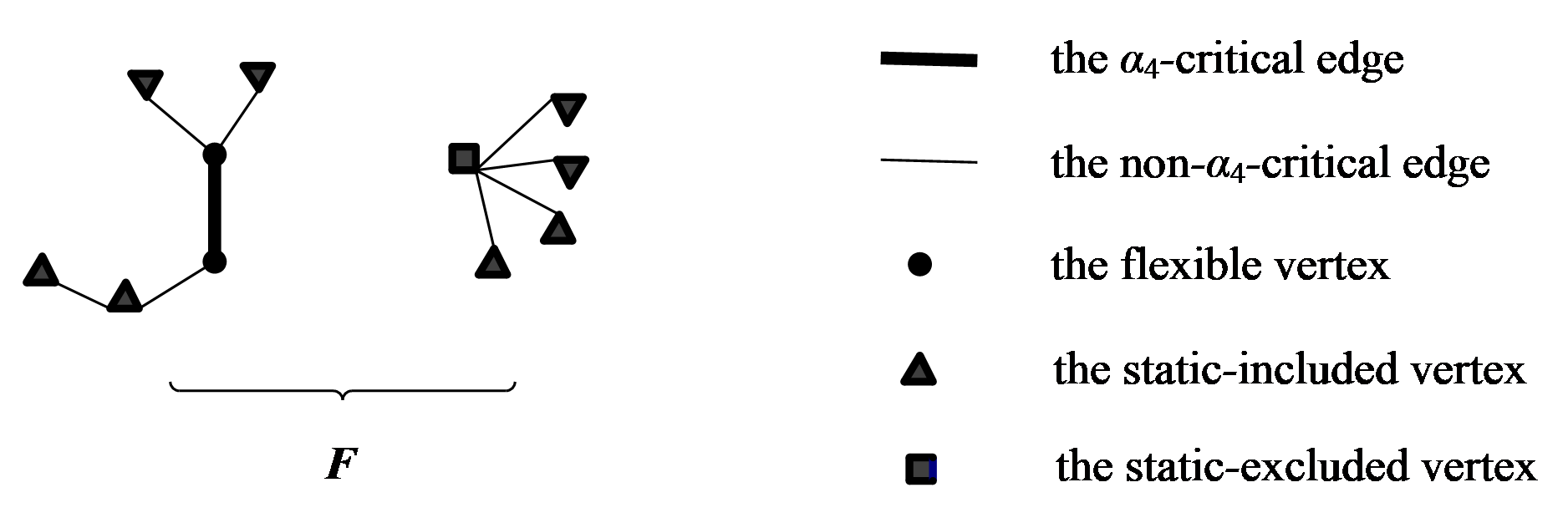}
    \caption{An example forest for $k=4$.}
    \label{fig:eg0} 
\end{figure}

Figure \ref{fig:eg0} shows an example forest $F$ for $k=4$, where $\mu_4(F)=2$ and $\alpha_4(F)=9$. For clarity, we use a thick solid line, a solid line, a solid circle, a solid triangle and a solid square to represent an $\alpha_k$-critical edge, a non-$\alpha_k$-critical edge, a flexible vertex, a static-included vertex and a static-excluded vertex, respectively, in all subsequent figures.

The concepts of flexible vertices and $\alpha_k$-critical edges are related. 

\begin{propositionNoParens}[\cite{3}]\label{pro:iff}
For a forest $F,$ a vertex $v$ is flexible in $F$ if and only if $v$ is incident with an $\alpha_k$-critical edge.
\end{propositionNoParens}

For a forest $F$, let $CE_k$ denote the set of all $\alpha_k$-critical edges in $F$, and let $F[CE_k]$ be the edge-induced subgraph formed by all $\alpha_k$-critical edges. By Proposition \ref{pro:iff}, $V(F[CE_k])=\mathcal{F}_F$. Therefore, $F[CE_k]$ is a subgraph of the vertex-induced subgraph $F[\mathcal{F}_F]$. Each connected component of $F[CE_k]$ (called an $\alpha_k$-\textit{critical component}) is a tree, and it is also called an $\alpha_k$-\textit{critical $t$-tree} if it has order $t$. We give the lower and upper bounds of $t$.

\begin{lemma}\label{lem:t}
Let $F$ be a forest and $H$ be an $\alpha_k$-critical $t$-tree of $F$. Then$,$ $2\leq t\leq k$.
\end{lemma}

\begin{proof}
Since $F[CE_k]$ is edge-induced and $H$ is a connected component of $F[CE_k]$, we obtain that $|V(H)|=t\geq2$. We next show that $t\leq k$. Suppose that $|V(H)|=t\geq k+1$. Let $M$ be an arbitrary maximum generalized $k$-mathcing of $F$, denoted by $\{T_1,\,T_2,\,\cdots,\,T_{\mu_{k}(F)}\}$. Clearly, $|V(T_i)|=k<|V(H)|$ for each $i$. Since $V(T_i)\cap V(T_j)=\varnothing$ for any $1\leq i<j\leq \mu_k(F)$ by the definition of $M$, we know that $E(H)\not\subseteq E(M)$. Therefore, there exists an edge $e\in E(H)$ such that $e\notin E(M)$. Note that $e$ is $\alpha_k$-critical by the definition of $H$; thus, $e\in E(M)$ according to Lemma \ref{lem:mu}, which is a contradiction. As a result, $2\leq t=|V(H)|\leq k$.
\end{proof}

Let $\Phi_k(F)$ be the total number of maximum generalized $k$-independent sets of $F$. The following theorem shows the relationship between the maximum generalized $k$-independent sets and $\alpha_k$-critical components in $F$.

\begin{theoremNoParens}[\cite{3}]\label{thm:V(H)-1}
Let $F$ be a forest and $H$ be an arbitrary $\alpha_k$-critical component of $F$. Then$,$ every maximum generalized $k$-independent set of $F$ contains exactly $|V(H)|-1$ vertices of $H$.
\end{theoremNoParens}

\begin{corollary}\label{cor:pi}
Let $F$ be a forest and $\mathcal{C}_F$ be the set of all $\alpha_k$-critical components of $F$. Then$,$ $\Phi_k(F)\leq \prod\limits_{H\in\mathcal{C}_F}|V(H)|$.
\end{corollary}

\begin{proof}
Denote $\mathcal{C}_F$ by $\{H_1,\,H_2,\,\cdots,\,H_{|\mathcal{C}_F|}\}$. Thus, all flexible vertices are included in $\bigcup\limits_{i=1}^{|\mathcal{C}_F|}|V(H_i)|$ by Proposition \ref{pro:iff}. Let $I$ be an arbitrary maximum generalized $k$-independent set of $F$. Clearly, $I$ is determined Completely by the flexible vertices in it, which means $I$ is determined by $\bigcup\limits_{i=1}^{|\mathcal{C}_F|}I\cap V(H_i)$. By Theorem \ref{thm:V(H)-1}, $|I\cap V(H_i)|=|V(H_i)|-1$ for each $i$. Thus, we obtain that $\Phi_k(F)\leq \prod\limits_{i=1}^{|\mathcal{C}_F|}|V(H_i)|$.
\end{proof}

The numbers of static-excluded vertices and $\alpha_k$-critical components in $F$ are related. 

\begin{lemmaNoParens}[\cite{3}]\label{lem:+}
Let $F$ be a forest and $\mathcal{C}_F$ be the set of all $\alpha_k$-critical components of $F$. Then, $|\mathcal{C}_F|+|\mathcal{N}_F|=\mu_k(F)$.
\end{lemmaNoParens}

\begin{remark}\label{rmk:+}
By the proof of Lemma \ref{lem:+}, for any maximum generalized $k$-matching $M$ of $F$, there exists a bijection between $M$ and $\mathcal{C}_F \cup \mathcal{N}_F$ such that each $k$-tree in $M$ either contains exactly one $\alpha_k$-critical component as a subgraph or contains exactly one static-excluded vertex in its vertex set.
\end{remark}

We introduce the following definitions to characterize the structure of an $\alpha_k$-critical component $H$ together with its adjacent component(s) in $F[\mathcal{A}_F]$.

\begin{definition}\label{AN[H]}
For a connected subgraph $T$ of a forest $F,$ its \textbf{$\mathcal{A}$-neighborhood} $N_{\mathcal{A}}(T)$ is the set of all connected components of $F[\mathcal{A}_F]$ that are adjacent to $T$. The \textbf{$\mathcal{A}$-closed neighborhood} of $T$, denoted by $N_{\mathcal{A}}[T]$, is the vertex-induced subgraph of $F$ on the vertex set $V(T)\cup V(N_{\mathcal{A}}(T))$.
\end{definition}

\begin{figure}[htbp]
    \centering 
    \includegraphics[width=1\textwidth]{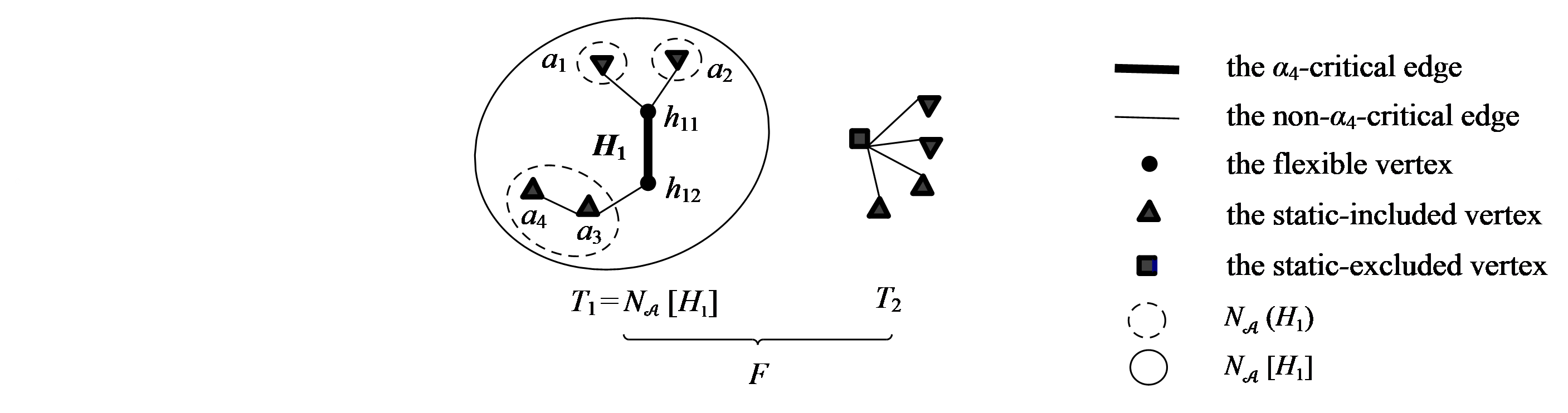}
    \caption{An example forest for $k=4$.} 
    \label{fig:eg0[]}
\end{figure}

In Figure \ref{fig:eg0[]}, $H_1$ is the $\alpha_4$-critical component. In fact, $N_{\mathcal{A}}(H_1)$ is the union of $N_{\mathcal{A}}(h_{11})$ and $N_{\mathcal{A}}(h_{12})$, which contains two vertices $\{a_1,\,a_2\}$ and a $K_2$ induced by $\{a_3, a_4\}$, respectively. Moreover, $N_{\mathcal{A}}[H_1]=T_1$. 

The following lemma was given by Li and Xu. We restate it with notations $N_{\mathcal{A}}(H)$ and $N_{\mathcal{A}}[H]$.

\begin{lemmaNoParens}[\cite{3}]\label{lem:a}
Let $F$ be a forest$,$ $H$ be an $\alpha_k$-critical $t$-tree of order $2\leq t\leq k-1,$ and $u$ be an arbitrary vertex of $H$. Then$,$

\noindent $(1)$ $|V(N_{\mathcal{A}}[H])|=|V(H)|+|V(N_{\mathcal{A}}(H))|\geq k\,;$

\noindent $(2)$ the largest component of $N_{\mathcal{A}}[H]-u$ has order at most $k-1\,;$

\noindent $(3)$ each leaf $v$ of $H$ satisfies $N_{\mathcal{A}}(v)\neq\varnothing$.
\end{lemmaNoParens}

Next, we give more properties of $N_{\mathcal{A}}(H)$ and $N_{\mathcal{A}}[H]$.

\begin{lemma}\label{lem:AN[H]}
Let $F$ be a forest$,$ and let $H$ be an $\alpha_k$-critical $t$-tree $(2\leq t\leq k)$ in $F.$ Then$,$

\noindent $(1)$ $N_{\mathcal{A}}[H]$ is a subtree of $F\,;$

\noindent $(2)$ $N_{\mathcal{A}}[H]=H$ if and only if $t=k$. And if $2\leq t\leq k-1$, then $N_{\mathcal{A}}[H]$ has order between $k+1$ and $2k-2$.
\end{lemma}

\begin{proof}
(1) By the definition, the induced subgraph $N_{\mathcal{A}}[H]$ is clearly connected.

(2) First, we need to show that $N_{\mathcal{A}}(H)=\varnothing$ if and only if $t=k$. Let $t=k$. Assume that $N_{\mathcal{A}}(H)\neq\varnothing$ and $A\in N_{\mathcal{A}}(H)$. Since $H$ is a tree with at least two leaves, we may choose a leaf $v\in V(H)$ that is not adjacent to $A$. Note that $v$ is a flexible vertex of $F$. Thus, there exists a maximum generalized $k$-independent set $I$ of $F$ such that $v\notin I$. By Theorem \ref{thm:V(H)-1}, $I\,\cap\,V(H)=V(H)\backslash\{v\}$. Moreover, $V(A)\subset I$ because all vertices of $A$ are static-included. Consequently, the subgraph induced by $V(A)\cup (V(H)\backslash\{v\})$ contains a tree of order at least $|V(H)|-1+|V(A)|\geq k-1+1=k$, contradicting the definition of $I$. Therefore, $N_{\mathcal{A}}(H)=\varnothing$. Conversely, if $N_{\mathcal{A}}(H)=\varnothing$ and we assume that $t\neq k$, then $2\leq t\leq k-1$ and $|V(N_{\mathcal{A}}[H])|=|V(H)|+|V(N_{\mathcal{A}}(H))|\leq k-1+0=k-1$, contradicting Lemma \ref{lem:a} (1). So $N_{\mathcal{A}}[H]=H$ if and only if $t=k$.

We now prove the second assertion. Let $2\leq t\leq k-1$. If $|V(N_{\mathcal{A}}[H])|\leq k$, by Lemma \ref{lem:a} (1), $N_{\mathcal{A}}[H]$ is a $k$-tree. Let $M$ be an arbitrary maximum generalized $k$-matching $M$ of $F$. By Lemma \ref{lem:+} and Remark \ref{rmk:+}, there exists a $k$-tree $T$ in $M$ that satisfies $V(H)\subset V(T)$ and $V(T)\subset V(N_{\mathcal{A}}[H])$. Hence, $T=N_{\mathcal{A}}[H]$. By the arbitrariness of $M$, every edge of $N_{\mathcal{A}}[H]$ is $\alpha_k$-critical. So $N_{\mathcal{A}}(H)=\varnothing$, and consequently, $N_{\mathcal{A}}[H]=H$ and $t=k$, which contradicts the statement that $t\leq k-1$.

If $|V(N_{\mathcal{A}}[H])|\geq 2k-1$, then let $u$ be a leaf of $H$ such that $|V(N_{\mathcal{A}}[u])|$ is largest among all the leaves of $H$. There is a component of $N_{\mathcal{A}}[H]-u$ whose vertex set is $V(N_{\mathcal{A}}[H])-V(N_{\mathcal{A}}[u])$, and it has order at most $k-1$ by Lemma \ref{lem:a} (2). So $|V(N_{\mathcal{A}}[H])|-|V(N_{\mathcal{A}}[u])|\leq k-1$. Hence, $2k-1\leq|V(N_{\mathcal{A}}[H])|\leq k+|V(N_{\mathcal{A}}[u])|-1$, which implies $|V(N_{\mathcal{A}}[u])|\geq k$. Let $I$ be an arbitrary maximum generalized $k$-independent set of $F$ containing $u$. Then, $I$ induces a tree $N_{\mathcal{A}}[u]$ on at least $k$ vertices, contradicting the definition of $I$. As a result, $k+1\leq|V(N_{\mathcal{A}}[H])|\leq 2k-2$.
\end{proof}

\begin{lemma}\label{lem:AN[H]super}
Let $F$ be a forest$,$ and let $H$ be an $\alpha_k$-critical component of $F.$ Let $k+q$ denote the order of $N_{\mathcal{A}}[H]$ $(0\leq q\leq k-2)$. Then$,$

\noindent $(1)$ the inequality $\min\limits_{u\text{ is a leaf of } H}|V(N_{\mathcal{A}}(u))|\geq q$ holds$\,;$ 

\noindent $(2)$ $|V(H)|\leq k-q$. If $q=0,$ then the equality always holds. If $1\leq q\leq k-2,$ then the equality holds if and only if $H$ is a $(k-q)$-path with ends $v_1$ and $v_2$ such that $N_{\mathcal{A}}(H)=N_{\mathcal{A}}(v_1)\cup N_{\mathcal{A}}(v_2)$ and $|N_{\mathcal{A}}(v_1)|=|N_{\mathcal{A}}(v_2)|=q$.
\end{lemma}

\begin{proof}
By Lemma \ref{lem:t} we know $2\leq|V(H)|\leq k$. 

(1) If $q=0$, then the inequality clearly holds. If $1\leq q\leq k-2$, and $\min\limits_{u\text{ is a leaf of } H}|V(N_{\mathcal{A}}(u))|\leq q-1$, then let $u$ be the corresponding leaf of $H$. The vertex set $V(N_{\mathcal{A}}[H])-V(N_{\mathcal{A}}[u])$ induces a connected subgraph of $N_{\mathcal{A}}[H]-u$, and this subtree has order $(k+q)-(1+|V(N_{\mathcal{A}}(u))|)\geq k$. However, by Lemma \ref{lem:a} (2) we know that each component of $N_{\mathcal{A}}[H]-u$ has order at most $k-1$, a contradiction. So $\min\limits_{u\text{ is a leaf of } H}|V(N_{\mathcal{A}}(u))|\geq q$. 

(2) If $q=0$, then by Lemma \ref{lem:AN[H]} (2), we know that $N_{\mathcal{A}}[H]=H$ has order $k$. So the equality $|V(H)|=k-q$ holds. If $1\leq q\leq k-2$, then by (1), we know that $\min\limits_{u\text{ is a leaf of } H}|V(N_{\mathcal{A}}(u))|\geq q$. Note that $V(N_{\mathcal{A}}[H])=V(H)\cup \bigcup\limits_{u\in V(H)}V(N_{\mathcal{A}}(u))$ and that $H$ has at least two leaves. Thus, $k+q=|V(N_{\mathcal{A}}[H])|\geq |V(H)|+2q$, leading to $|V(H)|\leq k-q$. If the equality holds, then $H$ has exactly two leaves $v_1$ and $v_2$ such that $N_{\mathcal{A}}(H)=N_{\mathcal{A}}(v_1)\cup N_{\mathcal{A}}(v_2)$ and $|N_{\mathcal{A}}(v_1)|=|N_{\mathcal{A}}(v_2)|=q$. Clearly, $H$ is a $(k-q)$-path.

Conversely, if $H$ is a $(k-q)$-path with ends $v_1$ and $v_2$ such that $N_{\mathcal{A}}(H)=N_{\mathcal{A}}(v_1)\cup N_{\mathcal{A}}(v_2)$ and $|N_{\mathcal{A}}(v_1)|=|N_{\mathcal{A}}(v_2)|=q$, then clearly $|V(H)|=k-q$.
\end{proof}

\begin{lemma}\label{lem:delta}
For a tree $T$ of order $n=k+q$ $(0\leq q\leq k-2),$ if $H$ is an $\alpha_k$-critical component of $T,$ then $T=N_{\mathcal{A}}[H]$ and $\Phi_k(T)=|V(H)|$.
\end{lemma}

\begin{proof}
Clearly, $\mu_k(T)=1$ for $n=k+q<2k$. Thus, $\alpha_k(T)=n-1$. Moreover, by Lemma \ref{lem:+}, $H$ is the only $\alpha_k$-component of $T$, and $\mathcal{N}_T=\varnothing$. Consequently, all vertices in $V(T)\backslash V(H)$ are static-included, leading to $T=N_{\mathcal{A}}[H]$. Let $u$ be an arbitrary vertex in $H$. By the definition of $H$, there exists a vertex $v\in V(H)$ such that $uv\in E(H)$. Clearly, $uv$ is $\alpha_k$-critical in $T$, so $\mu_k(T-uv)=0$. Thus, $T-u$ contains no subtrees of order $k$. As a result, $V(T)\backslash\{u\}$ is a maximum generalized $k$-independent set of $T$. So $\Phi_k(T)\geq|V(H)|$ by the arbitrariness of $u$. By Corollary \ref{cor:pi}, we obtain that $\Phi_k(T)\leq|V(H)|$. Thus, $\Phi_k(T)=|V(H)|$.
\end{proof}

Next we consider the relationships between different $\alpha_k$-critical components. We call two $\alpha_k$-critical components $H_1$ and $H_2$ \textbf{directly connected} if the subgraph of $F$ induced by $V(N_{\mathcal{A}}[H_1])\cup V(N_{\mathcal{A}}[H_2])$ is connected. Thus, there are two cases for directly connected $\alpha_k$-critical components $H_1$ and $H_2$: (1) there exists an edge $x_1x_2\in E(F)$ such that $x_1\in V(H_1)$ and $x_2\in V(H_2)$; (2) there exists a component $A$ of $F[\mathcal{A}_F]$ such that $A \in N_{\mathcal{A}}(H_1)\cap N_{\mathcal{A}}(H_2)$. In a forest $F$, there are no cycles. Thus, for two directly connected $\alpha_k$-critical components, exactly one of the two cases listed above holds, and not both.

\begin{figure}[htbp]
    \centering 
    \includegraphics[width=0.84\textwidth]{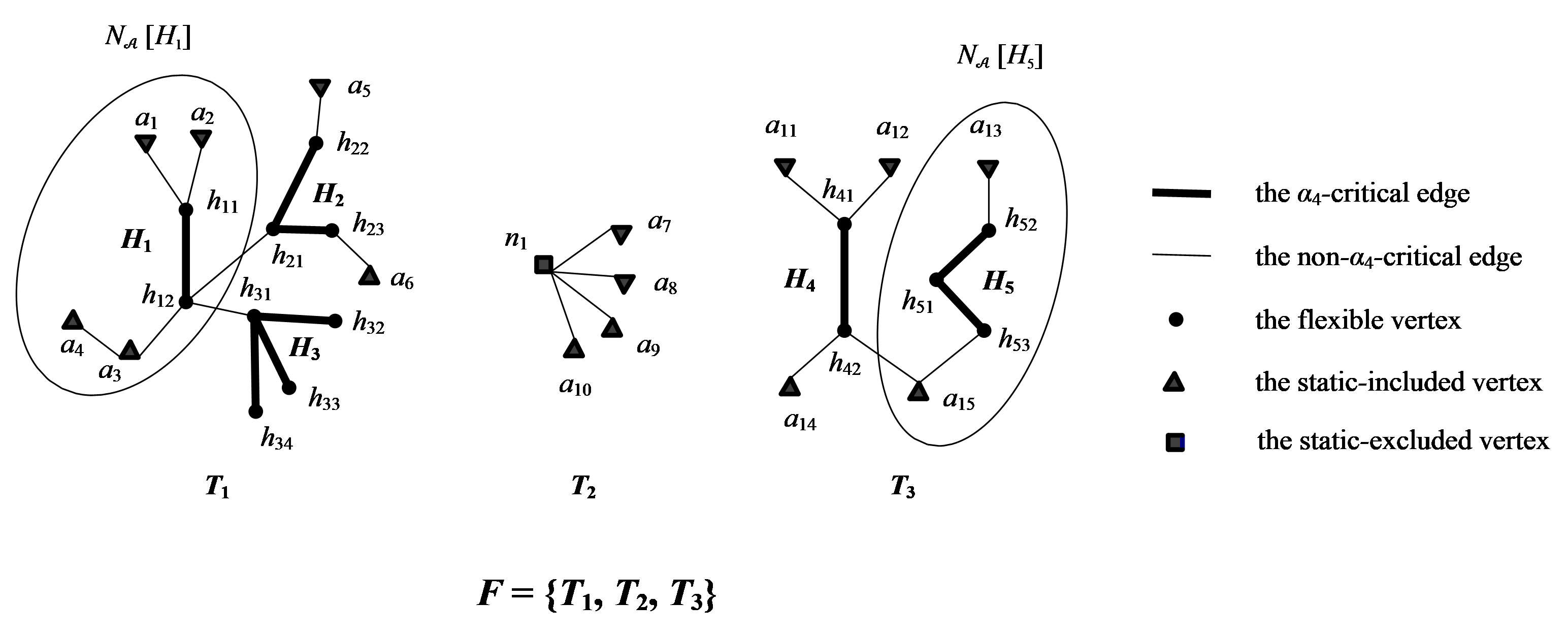}
    \caption{An example forest for $k=4$.}
    \label{fig:AN[H]}
\end{figure}

Figure \ref{fig:AN[H]} shows an example forest $F$ for $k=4$. A comprehensive examination of all the vertex types and $\alpha_4$-critical components in $F$ is provided in Appendix \ref{app:appendixA}. In this example, $H_1$ and $H_2$ are directly connected by edge $h_{12}h_{21}$. $H_4$ and $H_5$ are directly connected by $a_{15}$, which is a component of $F[\mathcal{A}_F]$. Moreover, $H_2$ and $H_3$ are not directly connected.

\section{Main Results}
\quad\,\,In this section, we establish the upper bound of $\Phi_k(F)$ for a forest $F$ and prove the main results, which are shown in Theorem \ref{thm:g-zzh} and Theorem \ref{thm:g-lsc}.

\subsection{Forests without Constraints}
\quad\,\,For clarity, we restate our result here.

\begin{reptheorem}[\ref{thm:g-zzh}]
Let $k\geq 2$ be an integer$,$ and let $F$ be a forest of order $n=pk+q$ $(n\geq1,\,0\leq q\leq k-1)$. Then$,$ $\Phi_k(F)\leq k^{p}$. The equality holds if and only if $F$ consists of $p$ trees of order $k$ and an arbitrary forest of order $q$.
\end{reptheorem}

\begin{proof}
Let $\mathcal{C}_F$ be the set of all $\alpha_k$-critical components in $F$, denoted by $\{H_1,\,H_2,\,\cdots,\,H_{|\mathcal{C}_F|}\}$. By Lemma \ref{lem:t}, $2\leq |V(H_i)|\leq k$ for each $1\leq i\leq |\mathcal{C}_F|$. By Corollary \ref{cor:pi}, we obtain

\begin{equation}\label{equ2}
\Phi_k(F)\leq \prod\limits_{i=1}^{|\mathcal{C}_F|}|V(H_i)|.
\end{equation}
By Lemma \ref{lem:+}, $|\mathcal{C}_F|=\mu_k(F)-|\mathcal{N}_F|\leq\mu_k(F)\leq \left\lfloor\displaystyle\frac{n}{k}\right\rfloor=p$. Thus, $\Phi_k(F)\leq k^{p}$.

If the equality holds, then $\mathcal{N}_F=\varnothing$ and $\mathcal{C}_F$ consists of $p$ $\alpha_k$-critical $k$-trees $H_1,\,H_2,\,\cdots,\,H_p$. Hence, $|\mathcal{F}_F|=pk$ and $|\mathcal{A}_F|=q$. By Lemma \ref{lem:AN[H]} (2), $N_{\mathcal{A}}(H_i)=\varnothing$ for $1\leq i\leq p$, which means that all flexible vertices are not adjacent to any static-included vertices in $F$. Thus, $F=F[\mathcal{A}_F]\cup F[\mathcal{F}_F]$. If $p=1$, then $F\cong H_1\cup F[\mathcal{A}_F]$ and the theorem is true. Next we may assume that $p\geq2$.

We now show that any two distinct members of $\{H_i :1\leq i\leq p\}$ are not directly connected. Without loss of generality, choose $H_1$ and $H_2$. Assume that $H_1$ and $H_2$ are directly connected by $x_1x_2$, where $x_1\in V(H_1)$ and $x_2\in V(H_2)$. Let $u_1\neq x_1$ and $u_2\neq x_2$ be two leaves of $H_1$ and $H_2$, respectively. Since the equality holds in (\ref{equ2}), there exists a maximum generalized $k$-independent set $I$ of $F$ such that $I\cap V(H_1)=V(H_1)\backslash\{u_1\}$ and $I\cap V(H_2)=V(H_2)\backslash\{u_2\}$ according to Theorem \ref{thm:V(H)-1}. However, $F[I]$ then contains a connected subgraph on at least $|V(H_1)|+|V(H_2)|-2=2k-2\geq k$ vertices, contradicting the definition of $I$. Therefore, $H_i$ and $H_j$ are not directly connected in $F$ for any $1\leq i<j\leq p$. As a result, $F$ is the union of $p$ trees $H_1,\,H_2,\,\cdots,\,H_p$ of order $k$ and a forest $F[\mathcal{A}_F]$ of order $q$.

Conversely, let $F$ be a forest consisting of $p$ trees $T_1,\,T_2,\,\cdots,\,T_p$ of order $k$ and an arbitrary forest $F'$ of order $q$. Since $q<k$, all vertices in $F'$ are static-included vertices of $F$. Clearly, for an arbitrary maximum generalized $k$-independent set $I$ of $F$, $I$ contains $k-1$ vertices in each $T_i$ ($1\leq i\leq p$). Note that $V(F)\backslash\bigcup\limits_{i=1}^{p}\{u_i:u_i\text{ is an arbitrary vertex of }T_i\}$ is a maximum generalized $k$-independent set. Thus, $\Phi_k(F)=k^p$.
\end{proof}

\subsection{Forests with Minimum Component Order $k$}

\quad\,\,For the rest of this section, we consider the forests whose each component has at least $k$ vertices.

\begin{lemma}\label{lem:twokindsofH}
Let $k_1\geq2$ and $T_1$ be an arbitrary $k_1$-tree. Let $k_2\geq3$ and $n=k_2+q,$ where $1\leq q\leq k_2-2$. Let $T_2$ be a tree of order $n,$ which is obtained from a $(k_2-q)$-path $P$ by attaching arbitrary forests $F_1$ and $F_2$ of order $q$ to its two ends$,$ respectively. Then$,$

\noindent $(1)$ $T_1$ is an $\alpha_{k_1}$-critical component itself.

\noindent $(2)$ $T_2$ has exactly one $\alpha_{k_2}$-critical component $P$.
\end{lemma}

\begin{proof}
Let $T_1$ be an arbitrary $k_1$-tree. Clearly, $\mu_{k_1}(T_1)=1$. In fact, each $e\in E(T_1)$ is $\alpha_{k_1}$-critical in $T_1$ because $T_1-e$ has no subtree of order $k$. So $T_1$ is an $\alpha_{k_1}$-critical $k_1$-tree itself.

Let $T_2$ be a tree of order $n$, which is obtained from a $(k_2-q)$-path $P$ by attaching arbitrary forests $F_1$ and $F_2$ of order $q$ to its two ends, respectively. Thus, $\mu_{k_2}(T_2)=1$ since $n<2k_2$. Note that $|V(P)|\geq2$, so $E(P)\neq\varnothing$. In fact, each $e\in E(P)$ is $\alpha_{k_2}$-critical in $T_2$ because $T_2-e$ has no subtree of order $k_2$. Moreover, each $e'\in E(T_2)\backslash E(P)$ is non-$\alpha_{k_2}$-critical in $T_2$, because $T_2-e'$ has at least one subtree of order $k_2$. So $P$ is the only $\alpha_{k_2}$-critical component of $T_2$.
\end{proof}

\begin{lemma}\label{thm:single}
Let $k\geq2$ and $T$ be a tree of order $n=k+q$ $(0\leq q\leq k-1)$. Then$,$ $\Phi_k(T)\leq k-q$. The equality holds if and only if 
\[T\cong
\begin{cases} 
\text{any } k\text{-tree}, & \text{if } q = 0\,; \\

\text{a tree obtained from a } (k-q)\text{-path by attaching arbitrary \quad} & ~\\

\text{  forests } F_1 \text{ and } F_2 \text{ of order } q \text{ to its two ends}, \text{ respectively}, & \text{if } 1 \leq q \leq k-2\,; \\

\text{any } (2k-1)\text{-tree}, & \text{if } q = k-1.
\end{cases}\]
\end{lemma}

\begin{proof}

Let $T$ be a tree of order $n=k+q$ $(0\leq q\leq k-1)$. Since $k\leq n=k+q<2k$, we obtain that $\mu_k(T)=1$ and $\alpha_k(T)=n-1$. Thus, if $I$ is an arbitrary maximum generalized $k$-independent set of $T$, then $I=V(T)\backslash\{u\}$ for some $u\in V(T)$. We next prove this lemma by following three cases. 

\noindent\textbf{Case 1.} $q=0$.

Let $T$ be an arbitrary $k$-tree. By Lemma \ref{lem:twokindsofH} (1), $T$ is an $\alpha_k$-critical component itself. Thus, $\Phi_k(T)=|V(T)|=k$ by Lemma \ref{lem:delta}. 

\noindent\textbf{Case 2.} $1\leq q\leq k-2$.

In this case, $k\geq3$. Since $\mu_k(T)=1$, we know that $|\mathcal{N}_T|=0$ or $1$ by Lemma \ref{lem:+}.

If $|\mathcal{N}_T|=1$, then $|\mathcal{C}_T|=0$, which means that $T$ contains no $\alpha_k$-critical components. Therefore, $\mathcal{F}_T=\varnothing$. It follows that $\Phi_k(T)=1<k-q$.

If $|\mathcal{N}_T|=0$, then $T$ contains exactly one $\alpha_k$-critical component $H$ and $\mathcal{N}_T=\varnothing$. Therefore, $T=N_{\mathcal{A}}[H]$. By Lemma \ref{lem:delta}, we know $\Phi_k(T)=|V(H)|=t$. By Lemma \ref{lem:AN[H]super} (2), we know $t\leq k-q$. Thus, we obtain that $\Phi_k(T)\leq k-q$. Moreover, if the equality holds, then $H$ is a $(k-q)$-path whose ends are $u$ and $v$ such that $|V(N_{\mathcal{A}}(u))|=|V(N_{\mathcal{A}}(v))|=q$ and $V(N_{\mathcal{A}}(u))\cup V(N_{\mathcal{A}}(v))=\mathcal{A}_T$, which means that $T$ is a tree obtained by attaching two (possibly distinct) forests $F_1$ and $F_2$ of order $q$ to the two endpoints of a $(k-q)$-path, respectively. 

Conversely, if $T$ is formed by attaching arbitrary forests $F_1$ and $F_2$ of order $q$ to each endpoint of a $(k-q)$-path $P$, then $P$ is the only $\alpha_k$-critical component of $T$ by Lemma \ref{lem:twokindsofH} (2). Thus, $\Phi_k(T)=|V(P)|=k-q$ by Lemma \ref{lem:delta}.

\noindent\textbf{Case 3.} $q= k-1$.

Let $T$ be an arbitrary tree of order $2k-1$. To prove the result, we next show that $T$ has no flexible vertices. Suppose that $y\in V(T)$ is a flexible vertex, then there exists another flexible vertex $x\in V(T)$ such that $xy$ is an $\alpha_k$-critical edge in $T$ by Proposition \ref{pro:iff}. Let $X$ and $Y$ be the components of $T-xy$ containing $x$ and $y$, respectively, and assume without loss of generality that $|V(X)|\geq|V(Y)|$. Thus, $|V(X)|\geq k$. Since $y$ is flexible in $T$, there exists a maximum generalized $k$-independent set $I$ of $T$ such that $y\notin I$. Note that $\alpha_k(T)=n-1$, we know $I=V(T)\backslash\{y\}$. However, $X$ is then a subtree of $T[I]$ of order at least $k$, a contradiction. Consequently, there are no flexible vertices in $T$. Therefore, $\Phi_k(T)=1=k-q$.
\end{proof}

\begin{lemma}\label{thm:k=2}
Let $F$ be a forest of order $2p+q$ $(p\geq1,0\leq q\leq 1)$ and each component of $F$ has order at least $2$. If $q=0,$ then $\Phi_2(F)\leq 2^p$ and equality holds if and only if $F\cong pK_2$. If $q=1,$ then $\Phi_2(F)\leq 2^{p-1}$ and equality holds if and only if $F$ is a forest constructed by attaching some $($possibly none$)$ of the $p-1$ copies of $K_2$ to the central vertex of a $3$-path $P$.
\end{lemma}

\begin{figure}[htbp]
    \centering 
    \includegraphics[width=0.75\textwidth]{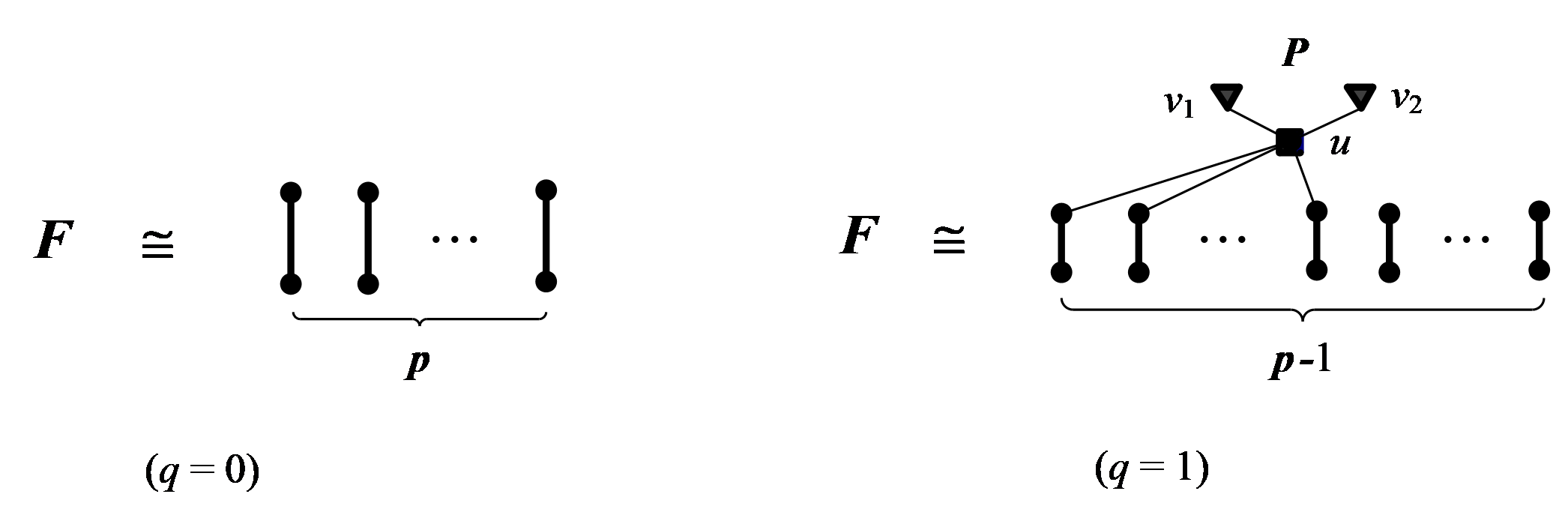}
    \caption{Extremal graphs for Lemma \ref{thm:k=2}.} 
    \label{fig:k=2} 
\end{figure}

\begin{proof}
Let $F$ be a forest of order $2p+q$ $(p\geq1,0\leq q\leq 1)$ and each component of $F$ has order at least $2$. For $q=0$, by Theorem \ref{thm:g-zzh}, $\Phi_2(F)\leq 2^p$ and equality holds if and only if $F\cong pK_2$.

For $q=1$, let $\mathcal{C}_F$ be the set of all $\alpha_2$-critical components in $F$. Let $H$ be an arbitrary $\alpha_2$-critical $t$-tree. By Lemma \ref{lem:t}, $t=2$. Consequently, each $\alpha_k$-critical component of $F$ is a $K_2$, which means that the number of $\alpha_2$-critical edges in $F$ is exactly $|\mathcal{C}_F|$. Thus, $|\mathcal{F}_F|=2|\mathcal{C}_F|$ by Proposition \ref{pro:iff}. Moreover, $N_{\mathcal{A}}[H]=H$ for each $H\in\mathcal{C}_F$ by Lemma \ref{lem:AN[H]}, which means that every static-included vertex and every flexible vertex are not adjacent in $F$. Note that in this lemma, any two static-included vertices are not adjacent in $F$ by the definitions of $\mathcal{A}_F$ and generalized $k$-independent set for $k=2$. Therefore, if $v\in\mathcal{A}_F$, then the neighborhood of $v$ in $F$ is composed of $d_F(v)$ static-excluded vertices. Moreover, $d_{F}(v)\geq1$ because each component of $F$ has order at least $2$.

By Lemma \ref{lem:+}, we obtain that $|\mathcal{C}_F|=\mu_2(F)-|\mathcal{N}_F|\leq\mu_2(F)\leq\lceil\frac{n}{2}\rceil = p$. If $|\mathcal{N}_F|=0$, then $|\mathcal{A}_F|=n-|\mathcal{N}_F|-|\mathcal{F}_F|=n-2|\mathcal{C}_F|\geq n-2p=1$. Let $v\in\mathcal{A}_F$; then, $|\mathcal{N}_F|\geq d_F(v)\geq1$, which is a contradiction. So $|\mathcal{N}_F|\geq1$ and $|\mathcal{C}_F|\leq p-1$. Thus, $\Phi_2(F)\leq 2^{p-1}$ by Corollary $\ref{cor:pi}$. If the equality holds, then $|\mathcal{C}_F|=p-1$ and $|\mathcal{N}_F|=1$. Moreover, $|\mathcal{F}_F|=2|\mathcal{C}_F|=2p-2$ and $|\mathcal{A}_F|=2$. Let $\mathcal{N}_F=\{u\}$ and $\mathcal{A}_F=\{v_1,v_2\}$. Thus, $uv_1\in E(F)$ and $uv_2\in E(F)$ by the discussion above. We next claim that arbitrary two elements of $\mathcal{C}_F$ are not joined. Assume that $x_1y_1$ and $x_2y_2$ are two $\alpha_2$-critical edges of $F$ such that $x_1x_2\in E(F)$. Since $\Phi_2(F)=2^{p-1}$, there exists a generalized 2-independent set $I$ of $F$ such that $x_1\in I$ and $x_2\in I$ by Theorem \ref{thm:V(H)-1}, contradicting the definition of $I$. As a result, $F$ is a forest constructed by attaching some (possibly none) of the $p-1$ copies of $K_2$ to the central vertex of a $3$-path $P$ (see Figure \ref{fig:k=2}).

Conversely, let $F$ be a forest constructed by attaching some (possibly none) of the $p-1$ copies of $K_2$ to the central vertex of a $3$-path $P$ (see Figure \ref{fig:k=2}). Then, for an arbitrary generalized $2$-independent set $I$ of $F$, $I$ contains at most one vertex from each of the $p-1$ copies of $K_2$ and at most $2$ vertices in $P$. Thus, each maximum generalized $2$-independent set of $F$ contains the two endpoints of $P$ and one arbitrary end of each $K_2$. So $\Phi_2(F)=2^{p-1}$.
\end{proof}

\begin{lemma}\label{lem:smallsuper}
Let $k\geq3$ be an integer. Let $T$ be a tree that satisfies $\mathcal{N}_T=\varnothing$ and contains at least two $\alpha_k$-critical components $\{H_1,\,H_2,\,\cdots\}$. If there exists an $H_i$ with $|V(N_{\mathcal{A}}[H_i])|:=k+q_i$ $(1\leq q_i\leq k-2)$ such that $H_i$ is directly connected to other $\alpha_k$-critical component$($s$)$ only through a single element $A\in N_{\mathcal{A}}(H_i),$ then $|V(A)|\leq q_i$.
\end{lemma}

\begin{figure}[htbp]
    \centering 
    \includegraphics[width=0.3\textwidth]{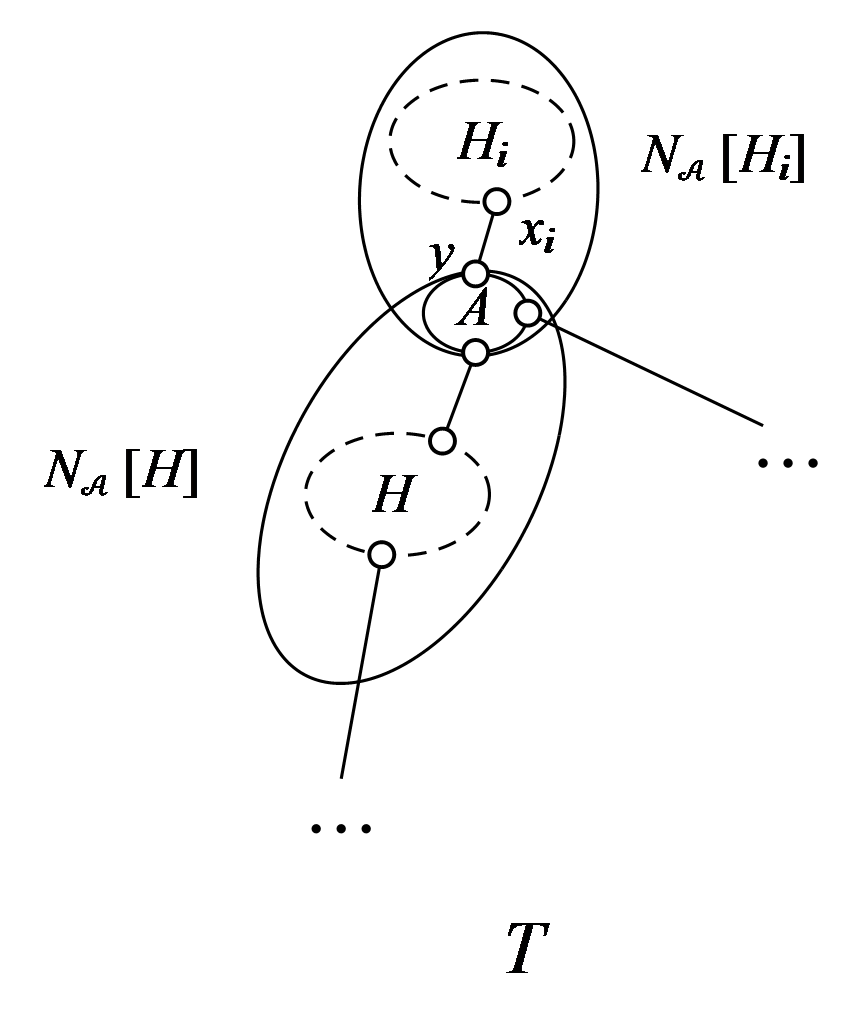}
    \vspace{-0.5cm}
    \caption{Illustration for Lemma \ref{lem:smallsuper}.} 
    \label{fig:smallsuper} 
\end{figure}

\begin{proof}
Assume that $|V(A)|\geq q_i+1$. Let $x_i\in V(H_i)$ and $y\in V(A)$ be the vertices of $T$ such that $A$ is attached to $H_i$ by edge $x_iy$ (see Figure \ref{fig:smallsuper}). Since $y$ is a static-included vertex of $T$, we know that $x_iy$ is non-$\alpha_k$-critical in $T$ by Proposition \ref{pro:iff}. Thus, $\mu_k(T-x_iy)=\mu_k(T)$. Let $T_{x_i}$ and $T_y$ denote the components of $T-x_iy$ containing $x_i$ and $y$, respectively. Thus, $T_{x_i}$ has order $|V(N_{\mathcal{A}}[H_i])|-|V(A)|=k+q_i-|V(A)|\leq k-1$. Thus, we know that $\mu_k(T_y)=\mu_k(T-x_iy)=\mu_k(T)$. Note that $T_y$ is a subgraph of $T-e$ for any $e\in E(H_i)$; we know $\mu_k(T-e)\geq\mu_k(T_y)=\mu_k(T)$, leading to $e$ being non-$\alpha_k$-critical in $T$. However, such an edge $e$ is $\alpha_k$-critical in $T$ by the definition of $H_i$, a contradiction. As a result, $|V(A)|\leq q_i$.
\end{proof}

\begin{lemma}\label{lem:super}
Let $k\geq3$ be an integer. Let $F$ be a forest whose each component has order at least $k$ such that $\mathcal{N}_F=\varnothing$. If there is a component $T$ of $F$ such that $T$ contains at least two $\alpha_k$-critical components of $F,$ then there exists a forest $F'$ whose each component has order at least $k$ such that $\mathcal{N}_{F'}=\varnothing$, $|V(F')|=|V(F)|,$ $\Phi_k(F')>\Phi_k(F),$ and $F'$ has the same number of $\alpha_k$-critical components as $F$. 
\end{lemma}

\begin{proof}
Without loss of generality, let $F=T$ be a tree of order at least $k$ such that $\mathcal{N}_T=\varnothing$. Let $\mathcal{C}_T=\{H_1,\,H_2,\,\cdots,\,H_p\}$ denote the $\alpha_k$-critical components of $T$, where $p\geq2$. Let $|V(N_{\mathcal{A}}[H_i])|:=k+q_i$ for each $1\leq i\leq p$, where $0\leq q_i\leq k-2$ by Lemma \ref{lem:AN[H]} (2). Thus, $|V(H_i)|\leq k-q_i$ for each $1\leq i\leq p$ by Lemma \ref{lem:AN[H]super} (2). Thus, by Corollary \ref{cor:pi} we obtain
\begin{equation}\label{equ2}
\Phi_k(T)\leq \prod\limits_{i=1}^{p}|V(H_i)|\leq (k-q_1)(k-q_2)\cdots(k-q_p).
\end{equation}

We claim that $\Phi_k(T)<(k-q_1)(k-q_2)\cdots(k-q_p).$ If not, then all the equalities in (\ref{equ2}) hold. Then, by Lemma \ref{lem:AN[H]super} (2) and Lemma \ref{lem:AN[H]} (2), we know that for each $i$, $N_{\mathcal{A}}[H_i]$ satisfies exactly one of the two conditions below: (a) $q_i=0$ and $N_{\mathcal{A}}[H_i]=H_i$ has order $k$; (b) $1\leq q_i\leq k-2$, $H_i$ is a path of order $k-q_i$ with endpoints $v_1^i$ and $v_2^i$ such that $N_{\mathcal{A}}(H_i)=N_{\mathcal{A}}(v_1^i)\cup N_{\mathcal{A}}(v_2^i)$ and $|N_{\mathcal{A}}(v_1^i)|=|N_{\mathcal{A}}(v_2^i)|=q_i$. Thus, if $q_i=0$ and $u_i$ is a leaf of $H_i$, then $V(N_{\mathcal{A}}[H_i])-V(N_{\mathcal{A}}[u_i])$ induces a subtree of order $k-1$. If $1\leq q_i\leq k-2$ and $u_i$ is a leaf of $H_i$, then $V(N_{\mathcal{A}}[H_i])-V(N_{\mathcal{A}}[u_i])$ still induces a subtree of order $k+q_i-(1+q_i)=k-1$. 

Since $T$ is connected and $|\mathcal{C}_T|=p\geq2$, there exist two $\alpha_k$-critical components, e.g., $H_1$ and $H_2$, that are directly connected in $T$. If $H_1$ and $H_2$ are directly connected by edge $x_1x_2$, where $x_1\in V(H_1)$ and $x_2\in V(H_2)$, then we choose two leaves $u_1\neq x_1$ and $u_2\neq x_2$ of $H_1$ and $H_2$, respectively. Thus, $u_1$ and $u_2$ are flexible in $T$ by the definitions of $H_1$ and $H_2$ according to Proposition \ref{pro:iff}. So there exists a maximum generalized $k$-independent set $I$ of $T$ such that $I\cap V(H_1)=V(H_1)\backslash\{u_1\}$ and $I\cap V(H_2)=V(H_2)\backslash\{u_2\}$ according to Theorem \ref{thm:V(H)-1} and the assumption that $\Phi_k(T)= \prod\limits_{i=1}^{p}|V(H_i)|$. Thus, $T[I]$ contains a subtree on the vertex set $(V(N_{\mathcal{A}}[H_1])-V(N_{\mathcal{A}}[u_1]))\cup (V(N_{\mathcal{A}}[H_2])-V(N_{\mathcal{A}}[u_2]))$. By the discussion above, we know that this subtree has order $(k-1)+(k-1)>k$, which contradicts the definition of $I$. If $H_1$ and $H_2$ are directly connected by $A\in N_{\mathcal{A}}(H_1)\cap N_{\mathcal{A}}(H_2)$, then $1\leq q_i\leq k-2$ for $i=1,\,2$. Moreover, $|V(A)|\leq q_1\leq k-2$. Similarly, $|V(A)|\leq q_2\leq k-2$. Let $x_1$ and $x_2$ be the vertices of $H_1$ and $H_2$ attached to $A$, respectively. Now, we choose two leaves $u_1\neq x_1$ and $u_2\neq x_2$ of $H_1$ and $H_2$, respectively. By similar discussion, there is a maximum generalized $k$-independent set $I$ of $T$ such that $T[I]$ contains a subtree on the vertex set $(V(N_{\mathcal{A}}[H_1])-V(N_{\mathcal{A}}[u_1]))\cup (V(N_{\mathcal{A}}[H_2])-V(N_{\mathcal{A}}[u_2]))$. We know that this subtree has order $(k-1)+(k-1)-|V(A)|\geq k$. Thus, $T[I]$ contains a subtree of order at least $k$, a contradiction. As a result, we obtain
\begin{equation}\label{equ3}
\Phi_k(T)<(k-q_1)(k-q_2)\cdots(k-q_p).
\end{equation}

Recall that $T$ has $p\geq2$ $\alpha_k$-critical components $\{H_1,\,H_2\cdots,\,H_p\}$. We define a set of subgraphs of $T$, that is, $S:=\{A: \text{ there exist } 1\leq i_1<i_2\leq p \text{ such that } A\in N_{\mathcal{A}}(H_{i_1})\cap N_{\mathcal{A}}(H_{i_2})\}=\{A_1,\,A_2,\,\cdots,\,A_{|S|}\}$. Let $T_C$ be a graph related to $T$ that has vertex set $V(T_C):=\{h_1,\,h_2,\,\cdots,\,h_p,$ $a_1,\,a_2,\,\cdots,\,a_{|S|}\}$ and edge set $E(T_C):=\{h_{i_1}h_{i_2}: H_{i_1} \text{ and } H_{i_2} \text{ are  directly connected by an edge} \}\cup\{h_{i}a_{j}: \text{ if } A_{j}\in N_{\mathcal{A}}(H_{i})\}$. Thus, $T_C$ is a tree, which shows all the location relationships of all the $\alpha_k$-critical components in $T$. We can observe that $d_{T_C}(a_{j})\geq2$ for each $1\leq j\leq |S|$. We give an example for $T$ and $T_C$ in Figure \ref{fig:TC}.

\begin{figure}[htbp]
    \centering 
    \includegraphics[width=0.5\textwidth]{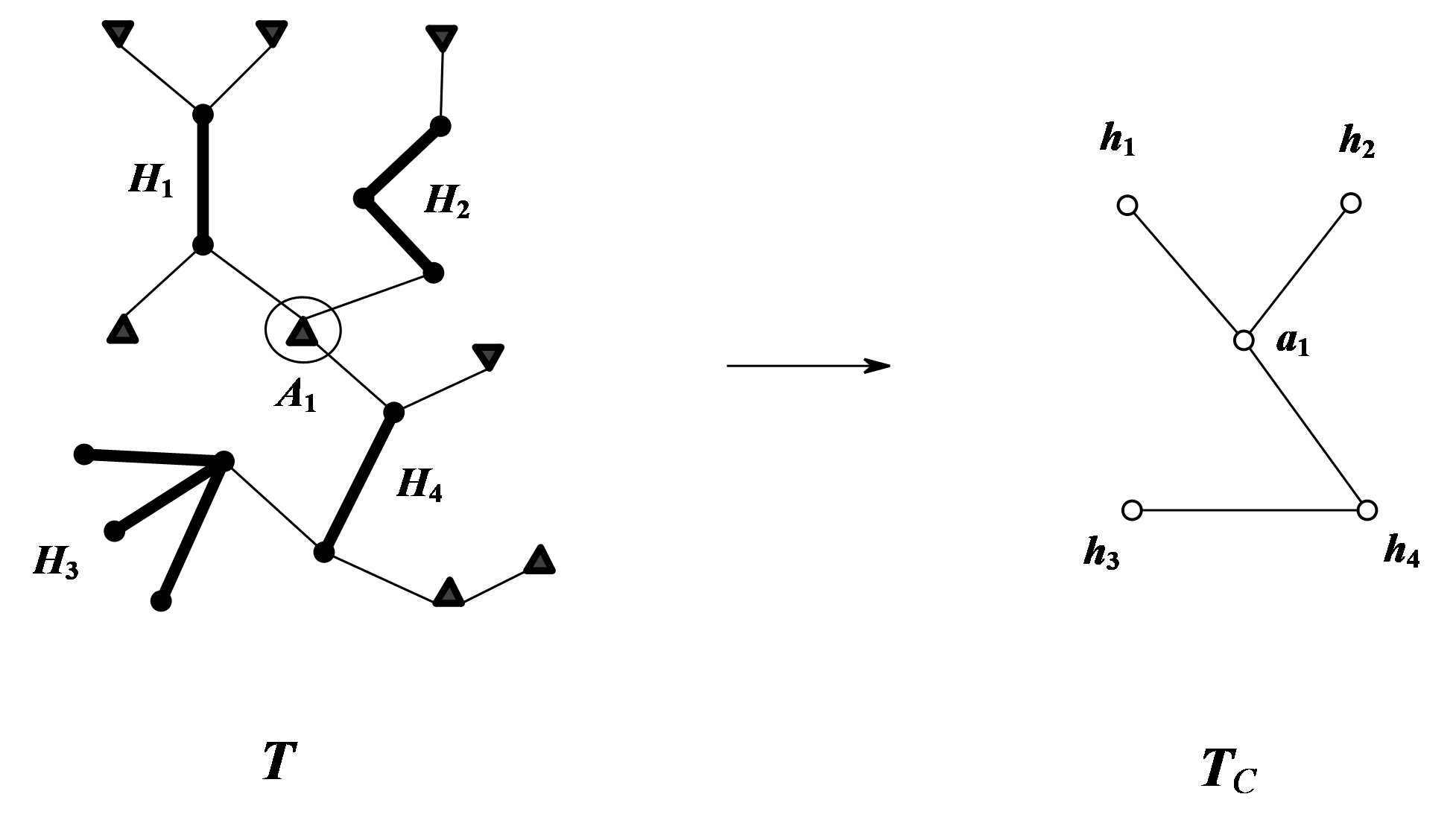}
    \caption{Illustration for an example of $T$ and $T_C$.} 
    \label{fig:TC} 
\end{figure}

Next we construct a forest $F'$ based on $T$, which satisfies the conditions that $\mathcal{N}_{F'}=\varnothing$, $|\mathcal{C}_{F'}|=|\mathcal{C}_{T}|=p$, and $|V(F')|=|V(T)|;$ additionally, each component of $F'$ has order at least $k$, and $\Phi_k(F')\geq(k-q_1)(k-q_2)\cdots(k-q_p)$.

Let $T_0:=T$, $G_0:=T_C$ and $F'_0=\varnothing$. There are total $p$ steps to construct the final forest $F'$, where the step $i$ ($1\leq i\leq p$) is defined as follows. 

\textbf{Step $1$}: Choose a leaf of $G_{0}$. Without loss of generality, let this leaf be $h_1$. If the neighbor of $h_1$ is $h_*$ in $G_{0}$, then let $x_1x_*$ be the edge of $T_{0}$ such that $x_1\in V(H_1)$ and $x_*\in V(H_*)$ (see Figure \ref{fig:step1}, where $i=1$). Then, let $G_{1}:=G_{0}-h_1$, let $T_{1}$ be the component of $T_{0}-x_1x_*$ which contains $x_*$, and let $F'_{1}:=F'_{0}\cup T'_{1}$, where $T'_{1}$ is a tree of order $|V(N_{\mathcal{A}}[H_1])|=k+q_1$ such that $\Phi_k(T'_1)=k-q_1$ (existing by Lemma \ref{thm:single}). If the neighbor of $h_1$ is $a_*$, then let $x_1\in V(H_1)$ and $y_*\in V(A_*)$ be the vertices of $T_{0}$ such that $H_1$ is attached to $A_*$ by edge $x_1y_*$ (see Figure \ref{fig:step1}, where $i=1$). Then, let $G_{1}=G_{0}-h_1$ if $d_{G_{0}}(a_*)\geq3$ or $G_{1}=G_{0}-h_1-a_*$ if $d_{G_{0}}(a_*)=2$. Let $T_{1}$ be the component of $T_{0}-x_1y_*$ which contains $y_*$, and let $F'_{1}:=F'_{0}\cup T'_{1}$, where $T'_{1}$ is a tree of order $|V(N_{\mathcal{A}}[H_1])|-|V(A_*)|=k+q_1-|V(A_*)|$ such that $\Phi_k(T'_1)=k-q_1+|V(A_*)|$ (existing by Lemma \ref{thm:single} and noting that $|V(A_*)|\leq q_1$ by Lemma \ref{lem:smallsuper}). Figure \ref{fig:realstep} shows an example for the definition of step 1.

\begin{figure}[htbp]
    \centering 
    \includegraphics[width=0.7\textwidth]{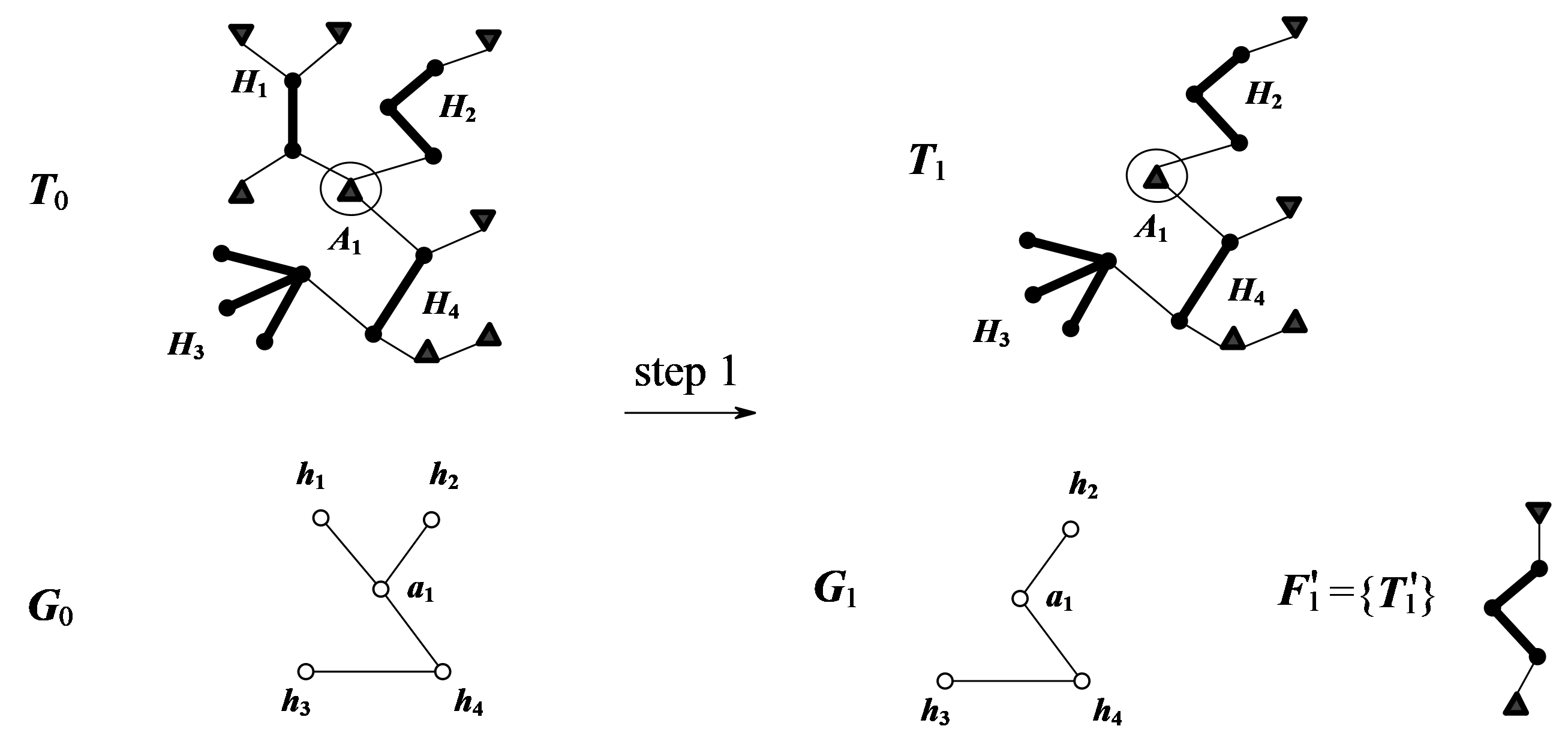}
    \caption{Illustration of an example of step 1.} 
    \label{fig:realstep} 
\end{figure}

\textbf{Step $i$} ($i\leq p-1$): Choose a leaf of $G_{i-1}$ and without loss of generality, let this leaf be $h_i$. Then, obtain $G_{i}$, $T_{i}$ and $F'_{i}$ based on $G_{i-1}$, $T_{i-1}$ and $F'_{i-1}$, respectively, by following the same procedures as in Step $1$, depending on the location of $h_i$ in $G_{i-1}$.

\begin{figure}[htbp]
    \centering 
    \includegraphics[width=0.9\textwidth]{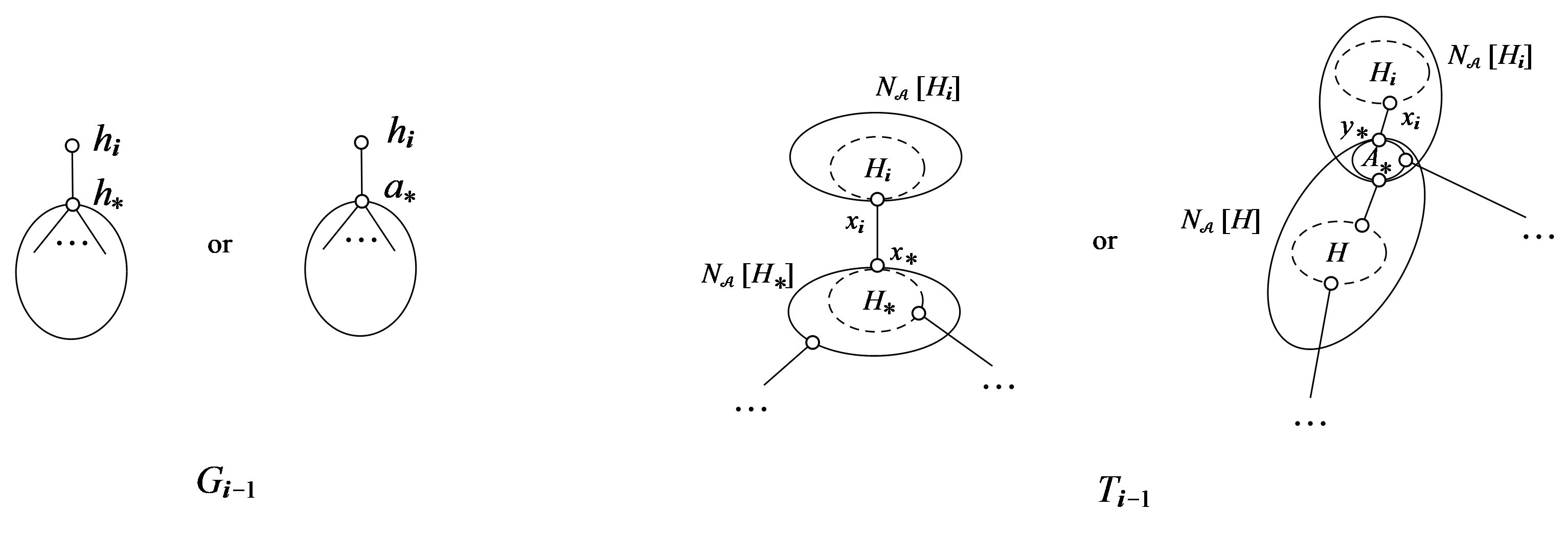}
    \caption{Illustration for the definition of Step $i$ ($i\leq p-1$).} 
    \label{fig:step1} 
\end{figure}

\textbf{Step $p$}: Now $|V(G_{p-1})|=1$, and without loss of generality, let this vertex be $h_p$. Let $G_{p}:=\varnothing$ and $T_{p}:=\varnothing$. Let $F'_{p}:=F'_{p-1}\cup T'_{p}$, where $T'_{p}$ is a tree of order $|V(N_{\mathcal{A}}[H_p])|=k+q_p$ such that $\Phi_k(T'_p)=k-q_p$ (existing by Lemma \ref{thm:single}).

Clearly, in each step $i$ ($1\leq i\leq p$), $|V(T'_i)|=|V(F'_i)|-|V(F'_{i-1})|=|V(T_{i-1})|-|V(T_{i})|$. Thus, we obtain $|V(T_{i-1})|+|V(F'_{i-1})|=|V(T_{i})|+|V(F'_i)|$, leading to $|V(T_{0})|+|V(F'_{0})|=|V(T_{p})|+|V(F'_{p})|$. Note that $|V(T_{p})|=0$ and $|V(F'_{0})|=0$; hence, we know that $|V(F'_{p})|=|V(T_{0})|=|V(T)|$. Note that $F'_{p}=\{T'_1,\,T'_2,\,\cdots,\,T'_p\}$ and $T'_i$ has order at least $k$ according to the process of each step $i$. Moreover, $\Phi_k(T'_i)$ is at least $k-q_i$. Thus, $F'_{p}$ is a forest whose each component has order at least $k$, and $\Phi_k(F'_{p})\geq(k-q_1)(k-q_2)\cdots(k-q_p)>\Phi_k(T)$ by (\ref{equ3}). 

Now, we show that $\mathcal{N}_{F'_{p}}=\varnothing$ and $|\mathcal{C}_{F'_{p}}|=|\mathcal{C}_{T}|=p$. By each step $1\leq i\leq p$, we know that each $T'_{i}$ is an extremal tree with order no more than $2k-2$ as described in Lemma \ref{thm:single}. Thus, $\mu_k(T'_i)=1$ and $T'_i$ has exactly one $\alpha_k$-critical component by Lemma \ref{lem:twokindsofH}. Therefore, $\mathcal{N}_{T'_{i}}=\varnothing$ by Lemma \ref{lem:+}. As a result, we know $\mathcal{N}_{F'_{p}}=\varnothing$ and $|\mathcal{C}_{F'_{p}}|=p$. Thus, we have constructed the forest $F':=F'_{p}$, completing the proof.
\end{proof}

Now, we prove another main result Theorem \ref{thm:g-lsc}. We restate it here for clarity.

\begin{reptheorem}[\ref{thm:g-lsc}]
Let $k\geq2$ be an integer$,$ and let $F$ be a forest of order $n=pk+q$ $(p\geq1,\,0\leq q\leq k-1),$ where each component of $F$ has order at least $k$. Then$,$ $\Phi_k(F)\leq (k-\lfloor\frac{q}{p}\rfloor)^{p-r}\cdot(k-\lceil\frac{q}{p}\rceil)^{r},$ where $0\leq r\leq p-1$ and $r\equiv q \pmod{p}$. Furthermore$,$ all the extremal graphs are described as follows$\,:$

\noindent $(1)$ If $k=2$ and $q=0,$ then the equality holds $\iff F\cong pK_2$.

\noindent $(2)$ If $k=2$ and $q=1,$ then the equality holds $\iff F$ is a forest constructed by attaching some $($possibly none$)$ of the $p-1$ copies of $K_2$ to the central vertex of a $3$-path $P$.

\noindent $(3)$ If $k\geq3$ and $p\geq q,$ then the equality holds $\iff F\cong (\bigcup\limits_{i=1}^{p-q} T_i)\cup qP_{k+1},$ where each $T_i$ is a $k$-tree.

\noindent $(4)$ If $k\geq3$ and $p<q$ with $p\neq1$ or $q\neq k-1,$ then the equality holds $\iff F\cong \bigcup\limits_{i=1}^{p} T_i,$ where each $T_i$ is a $(k+t_i)$-tree obtained from a $(k-t_i)$-path by attaching arbitrary forests $F_{i}^{1}$ and $F_{i}^{2}$ of order $t_i$ to its two ends$,$ respectively. Here$,$ $t_i=\lfloor\frac{q}{p}\rfloor$ for $1\leq i\leq p-r$ and $\lceil\frac{q}{p}\rceil$ for $p-r< i\leq p$.

\noindent $(5)$ If $k\geq3,$ $p=1$ and $q=k-1,$ then the equality holds $\iff F$ is an arbitrary $(2k-1)$-tree.
\end{reptheorem}

\begin{figure}[htbp]
    \centering 
    \includegraphics[width=1\textwidth]{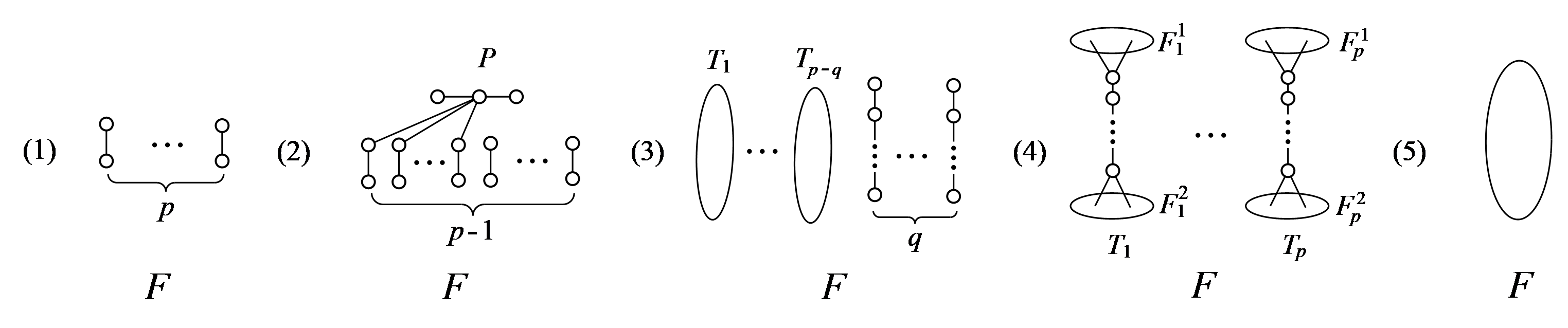}
    \caption{Extremal graphs for Theorem \ref{thm:g-lsc}.} 
\end{figure}

\begin{proof}
First, if $k=2$, then the theorem holds according to Lemma \ref{thm:k=2}. So we next consider $k\geq3$. If $p=1$, then the theorem again holds according to Lemma \ref{thm:single}. We now consider the case when $k\geq3$ and $p\geq2$. 

Let $\mathfrak{F}(n)$ denote the set of all forests of order $n=pk+q$ ($p\geq 2,\,0\leq q\leq k-1$), where each component has order at least $k$. Let $F\in \mathfrak{F}(n)$ be a forest. By Lemma \ref{lem:+}, we know that $|\mathcal{C}_F|=\mu_k(F)-|\mathcal{N}_F|\leq \mu_k(F)\leq \Big\lfloor\displaystyle\frac{n}{k}\Big\rfloor\leq p$. 

If $|\mathcal{C}_F|\leq p-1$, then by Corollary \ref{cor:pi}, we obtain $\Phi_k(F)\leq\prod\limits_{H\in\mathcal{C}_F}|V(H)|\leq k^{p-1}$ since $|V(H)|\leq k$ by Lemma \ref{lem:t}. When $q=0$, $k^{p-1}<k^p=(k-\lfloor\frac{q}{p}\rfloor)^{p-r}\cdot(k-\lceil\frac{q}{p}\rceil)^{r}$. When $q=1$, $k^{p-1}<k^{p-1}\cdot(k-1)=(k-\lfloor\frac{q}{p}\rfloor)^{p-r}\cdot(k-\lceil\frac{q}{p}\rceil)^{r}$. When $2\leq q \leq k-1$, $k^{p-1}\leq k^{p-1}\cdot(k-q)<(k-\lfloor\frac{q}{p}\rfloor)^{p-r}\cdot(k-\lceil\frac{q}{p}\rceil)^{r}$. The last inequality holds by applying the inequality $(k-a)(k-b)\geq k(k-a-b)$ for $a,b\geq0$, which is strict if $a>0$ and $b>0$. Consequently, $\Phi_k(F)\leq\prod\limits_{H\in\mathcal{C}_F}|V(H)|\leq k^{p-1}<(k-\lfloor\frac{q}{p}\rfloor)^{p-r}\cdot(k-\lceil\frac{q}{p}\rceil)^{r}$. 

So next we only need to consider the forest $F\in\mathfrak{F}(n)$ such that $|\mathcal{C}_F|=\mu_k(F)=p$ and $\mathcal{N}_{F}=\varnothing$. Choose an arbitrary forest $\widetilde{F}$ such that $\Phi_k(\widetilde{F})$ is the largest among such forests mentioned above. Let $\mathcal{C}_{\widetilde{F}}=\{H_1,\,H_2,\,\cdots,\,H_p\}$. Since each component of $\widetilde{F}$ has order at least $k$, we know that there are no components of $\widetilde{F}$ whose vertices are all static-included. Thus, each component of $\widetilde{F}$ contains at least one element of $\{H_1,\,H_2,\,\cdots,\,H_p\}$ as its subgraph. By Lemma \ref{lem:super}, we know that $\widetilde{F}$ has no component that contains at least two $\alpha_k$-critical components. Thus, $\widetilde{F}$ has exactly $p$ components $T_1,\,T_2,\,\cdots,\,T_p$ and $T_i=N_{\mathcal{A}}[H_i]$ for each $1\leq i\leq p$.

Let $T_i=N_{\mathcal{A}}[H_i]$ have order $k+q_i$. Thus, $0\leq q_i\leq k-2$ by Lemma \ref{lem:AN[H]} (2). Moreover, $\sum\limits_{i=1}^{p}q_i=n-pk=q$. Thus, by the choice of $\widetilde{F}$, we get $\Phi_k(T_i)= k-q_i$ according to Lemma \ref{thm:single}. Therefore, the following equality holds.
\begin{equation}\label{equ4}
\Phi_k(\widetilde{F})=(k-q_1)(k-q_2)\cdots(k-q_p).
\end{equation}
Moreover, $T_i$ is a $k$-tree (for $q_i=0$) or a $(k+q_i)$-tree obtained by attaching two (possibly distinct) forests $F_i^1$ and $F_i^2$ of order $q_i$ to the two endpoints of a $(k-q_i)$-path, respectively (for $1\leq q_i\leq k-2$).

In fact, $|q_i-q_j|\leq 1$ for any $1\leq i<j\leq p$. If not, let $q_1-q_2\geq2$ without loss of generality. By Lemma \ref{thm:single}, there exist two trees $T'_1$ and $T'_2$ such that $V(T'_1)=k+q_1-1$, $\Phi_k(T'_1)=k-q_1+1$, $V(T'_2)=k+q_2+1$ and $\Phi_k(T'_2)=k-q_2-1$. Moreover, $T'_1$ contains exactly one $\alpha_k$-critical component by Lemma \ref{lem:twokindsofH}, as does $T'_2$. Let $F':=(\widetilde{F}\backslash\{T_1,T_2\})\cup\{T'_1,T'_2\}$. Thus, $F'\in\mathfrak{F}(n)$, $|\mathcal{C}_{F'}|=|\mathcal{C}_{\widetilde{F}}|=p$ and $\mathcal{N}_{F'}=\varnothing$. Note that $\Phi_k(F')=(k-q_1+1)(k-q_2-1)\cdots(k-q_p)>(k-q_1)(k-q_2)\cdots(k-q_p)=\Phi_k(\widetilde{F})$, contradicting the choice of $\widetilde{F}$.

Note that $\{\lfloor\frac{q}{p}\rfloor,\,\cdots,\,\lfloor\frac{q}{p}\rfloor,\,\lceil\frac{q}{p}\rceil,\,\cdots,\,\lceil\frac{q}{p}\rceil\}$ is the unique sequence consisting of $p$ nonnegetive integers $\{q_1,\,q_2,\,\cdots,\,q_p\}$ such that $\sum\limits_{i=1}^{p}q_i=q$ and $|q_i-q_j|\leq1$ for any $1\leq i<j\leq p$. So $\Phi_k(\widetilde{F})=(k-\lfloor\frac{q}{p}\rfloor)^{p-r}\cdot(k-\lceil\frac{q}{p}\rceil)^{r}$, where $0\leq r\leq p-1$ and $r\equiv q \pmod{p}$. As a result, we know that $\Phi_k(F)\leq (k-\lfloor\frac{q}{p}\rfloor)^{p-r}\cdot(k-\lceil\frac{q}{p}\rceil)^{r}$ for an arbitrary $F\in\mathfrak{F}(n)$. By the arbitrariness of $\widetilde{F}$, if the equality holds for $p\geq q$, then $F$ contains $(p-q)$ trees of order $k$ and $q$ paths of order $k+1$. If the equality holds for $p<q$, then $F$ contains $p$ trees, and each tree has order $k+\lfloor\frac{q}{p}\rfloor$ (or $k+\lceil\frac{q}{p}\rceil$), which is obtained from a path of order $k-\lfloor\frac{q}{p}\rfloor$ (or $k-\lceil\frac{q}{p}\rceil$) by attaching arbitrary forests of order $\lfloor\frac{q}{p}\rfloor$ (or $\lceil\frac{q}{p}\rceil$) to its two ends, respectively.

Conversely, for the case when $p\geq q$, let $F\cong (\bigcup\limits_{i=1}^{p-q} T_i)\cup qP_{k+1},$ where each $T_i$ is a $k$-tree. Thus, by Lemma \ref{thm:single} we get $\Phi_k(F)=(k-0)^{p-q}(k-1)^q= (k-\lfloor\frac{q}{p}\rfloor)^{p-r}\cdot(k-\lceil\frac{q}{p}\rceil)^{r},$ where $0\leq r\leq p-1$ and $r\equiv q \pmod{p}$. For the case when $p<q$, let $F\cong \bigcup\limits_{i=1}^{p} T_i,$ where each $T_i$ is a $(k+t_i)$-tree obtained from a $(k-t_i)$-path by attaching arbitrary forests $F_{i}^{1}$ and $F_{i}^{2}$ of order $t_i$ to its two ends$,$ respectively. Here$,$ $t_i=\lfloor\frac{q}{p}\rfloor$ for $1\leq i\leq p-r$ and $\lceil\frac{q}{p}\rceil$ for $p-r< i\leq p$. Thus, $\Phi_k(T_i)=k-t_i$ by Lemma \ref{thm:single}. Consequently, $\Phi_k(F)= (k-\lfloor\frac{q}{p}\rfloor)^{p-r}\cdot(k-\lceil\frac{q}{p}\rceil)^{r},$ where $0\leq r\leq p-1$ and $r\equiv q \pmod{p}$. This completes the proof.
\end{proof}

\section{Conclusion and Remark}

\quad\,\,Our structural theorems refine the theory concerning the maximum generalized $k$-independent sets in trees proposed by Li and Xu \cite{3}, providing a more complete characterization of their properties. Building upon these properties, we conduct further analyses and computations.

For a general integer $k\geq2$, we establish the upper bound of $\Phi_k(F)$ of an arbitrary forest $F$ in Theorem \ref{thm:g-zzh}, and we characterize all the extremal forests. We also study an arbitrary forest $F$ such that each component of $F$ has order at least $k$, and establish the upper bound of $\Phi_k(F)$ in Theorem \ref{thm:g-lsc}. We again characterize all the corresponding extremal forests. These two results extend some works by Jou and Chang \cite{7} and Sun and Li \cite{8} from special integers to general integers, respectively.

\newpage

\begin{appendices}
\section{}\label{app:appendixA}
\begin{figure}[htbp]
    \centering 
    \includegraphics[width=1\textwidth]{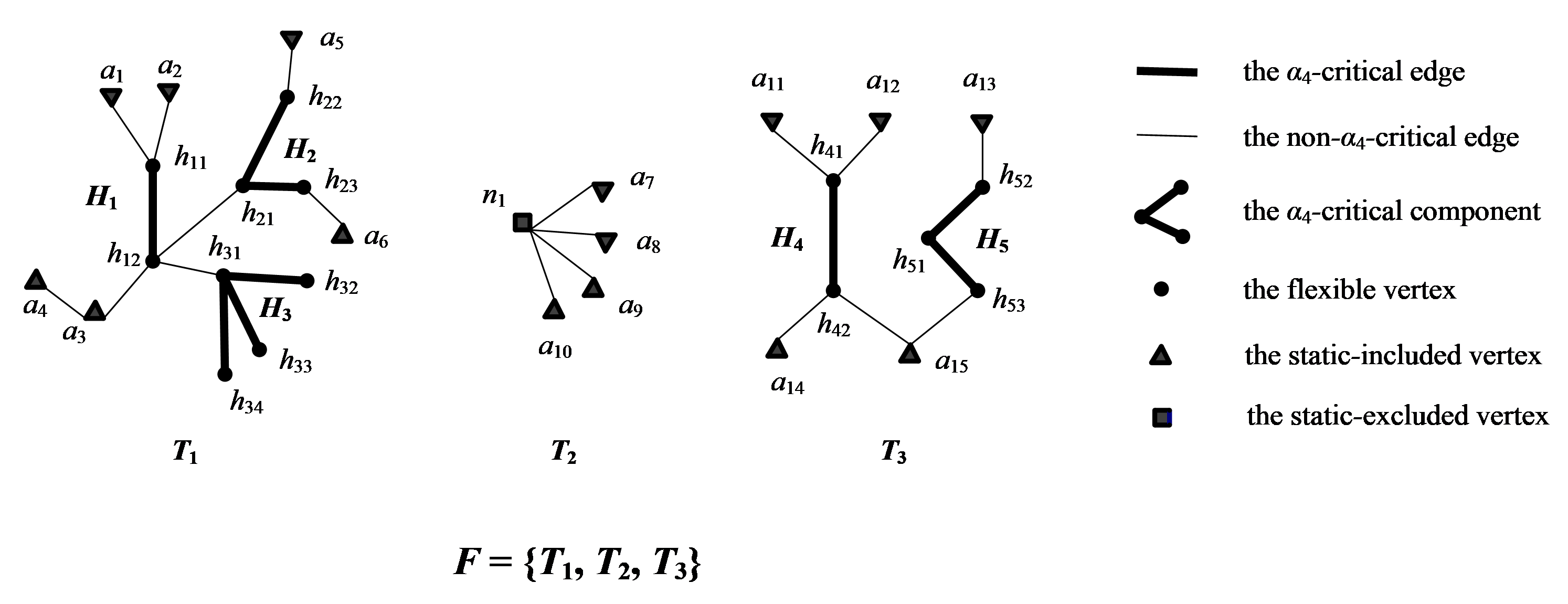}
    \caption{An example forest for $k=4$} 
    \label{fig:3vertices-eg2} 
\end{figure}

Let $F=\{T_1,\,T_2,\,T_3\}$ be the forest depicted in Figure \ref{fig:3vertices-eg2}, and set $k=4$. We now exemplify the three vertex types and all $\alpha_k$-critical components in $F$. Since $|V(T_1)|=15$, $|V(T_2)|=5$ and $|V(T_3)|=10$, we have that $\mu_4(T_1)\leq3$, $\mu_4(T_2)\leq1$ and $\mu_4(T_3)\leq2$; therefore, $\mu_4(F)=\sum\limits_{i=1}^3\mu_4(T_i)\leq6$. Conversely, one can easily construct a generalized $4$-matching of $F$ containing six vertex-disjoint 4-trees, so $\mu_4(F)=6$. 

Each edge shown as a thick solid line in Figure \ref{fig:3vertices-eg2} is $\alpha_4$-critical. For illustration purposes, consider the edge $h_{21}h_{22}$. Deleting it reduces the maximum possible number of disjoint $4$-trees in $T_1-h_{21}h_{22}$ to at most two, because any generalized $4$-matching $M$ of $T_1-h_{21}h_{22}$ meets at most one vertex of $\{a_1,\,a_4,\,a_6\}$, leading to $\mu_4(T_1-h_{21}h_{22})\leq \lfloor\frac{|V(T_1)\backslash\{a_5,\,h_{22}\}|-2}{4}\rfloor=2$. Hence $\mu_4(F-h_{21}h_{22})=\mu_4(T_1-h_{21}h_{22})+\mu_4(T_2)+\mu_4(T_3)\leq2+1+2=5<6=\mu_4(F)$. 

Each edge $e$ shown as a solid line (for example, $a_1h_{11}$, $h_{12}h_{21}$, and $n_1a_7$) in Figure \ref{fig:3vertices-eg2} is a non-$\alpha_4$-critical edge of $F$ because $\mu_4(F-e)=6=\mu_4(F)$. The forest $F$ therefore possesses exactly five $\alpha_k$-critical components $H_1,\,H_2,\,\cdots,H_5$. By Theorem \ref{thm:V(H)-1}, every maximum generalized $4$-independent set of $F$ contains $\alpha_4(F)-\sum\limits_{i=1}^5|V(H_i)-1|=|V(F)|-\mu_4(F)-\sum\limits_{i=1}^5|V(H_i)-1|=15$ vertices from the set $\{a_1,\,a_2,\,\cdots,\,a_{15},\,n_1\}$. A direct inspection shows that $n_1$ is a static-excluded vertex of $T_2$ (and hence of $F$). Thus, all three vertex types together with the $\alpha_4$-critical components in $F$ have been verified.
\end{appendices}

\end{document}